\documentclass[11pt]{article}
\usepackage[top=1in, bottom=1in, left=1in, right=1in]{geometry}
\usepackage{tikz}
\usetikzlibrary{shapes,arrows,matrix,patterns,positioning}
\usepackage{color}
\usepackage{graphicx}
\usepackage{algorithmic}
\usepackage{amssymb}
\usepackage{epstopdf}
\usepackage[ruled]{algorithm2e}
\usepackage{float}
\usepackage{amsmath,amsthm,bm,color,epsfig,enumerate,caption}
\usepackage{hyperref}
\usepackage{booktabs}
\usepackage{multirow}
\usepackage{array}

\newtheorem{theorem}{Theorem}[section]

\newtheorem{algo}[theorem]{Algorithm}

\renewcommand{\appendix}[1]{
\section*{Appendix: #1}
}
\newcommand{\norm}[1]{\left\lVert#1\right\rVert}
\renewcommand{\O}{O}

\newcommand{\bbC}{\mathbb{C}}
\newcommand{\bbR}{\mathbb{R}}

\makeatletter
\newcommand*{\extendadd}{
  \mathbin{
    \mathpalette\extend@add{}
  }
}
\newcommand*{\extend@add}[2]{
  \ooalign{
    $\m@th#1\leftrightarrow$%
    \vphantom{$\m@th#1\updownarrow$}
    \cr
    \hfil$\m@th#1\updownarrow$\hfil
  }
}
\makeatother

\begin{document}

\title{Butterfly Factorization}

\author{Yingzhou Li$^\sharp$,
    Haizhao Yang$^\dagger$,
    Eileen R. Martin$^\sharp$,
    Kenneth L. Ho$^\dagger$,
    Lexing Ying$^{\dagger\sharp}$
  \vspace{0.1in}\\
  $\dagger$ Department of Mathematics, Stanford University\\
  $\sharp$ ICME, Stanford University
}

\maketitle

\begin{abstract}
The paper introduces the butterfly factorization as a data-sparse
approximation for the matrices that satisfy a complementary low-rank
property. The factorization can be constructed efficiently if either
fast algorithms for applying the matrix and its adjoint are available
or the entries of the matrix can be sampled individually. For an
$N\times N$ matrix, the resulting factorization is a product of
$O(\log N)$ sparse matrices, each with $O(N)$ non-zero entries. Hence,
it can be applied rapidly in $O(N\log N)$ operations. Numerical
results are provided to demonstrate the effectiveness of the butterfly
factorization and its construction algorithms.
\end{abstract}



\section{Introduction}
\label{sec:intro}

One of the key problems in scientific computing is the rapid
evaluation of dense matrix-vector multiplication. Given a matrix $K
\in \bbC^{N\times N}$ and a vector $g\in\bbC^N$, the direct
computation of the vector $u=Kg\in \bbC^N$ takes $O(N^2)$ operations
since each entry of $K$ contributes to the result. This type of dense
multiplication problem appears widely in special function transforms,
integral transforms and equations, data fitting and smoothing,
etc. Since the number of unknowns $N$ is typically quite large in
these applications, a lot of work has been devoted to performing this
computation more efficiently without sacrificing the accuracy. Such a
reduction in computational complexity depends highly on the algebraic
and numerical properties of the matrix $K$. For certain types of
matrices $K$, such as the Fourier matrix, numerically low-rank
matrices, hierarchically semi-separable ({HSS}) matrices \cite{HSS1},
and hierarchical matrices \cite{HSS2,HMat}, there exist fast algorithms for
computing $Kg$ accurately in $O(N\log N)$ or even $O(N)$ operations.

\subsection{Complementary low-rank matrices and butterfly algorithm}

Recent work in this area has identified yet another class of matrices
for which fast $O(N\log N)$ application algorithms are available.
These matrices satisfy a special kind of complementary low-rank
property. For such a matrix, the rows are typically indexed by a set
of points, say $X$, and the columns by another set of points, say
$\Omega$. Both $X$ and $\Omega$ are often point sets in $\bbR^d$ for
some dimension $d$. Associated with $X$ and $\Omega$ are two trees
$T_X$ and $T_\Omega$, respectively and both trees are assumed to have
the same depth $L=O(\log N)$, with the top level being level $0$ and
the bottom one being level $L$. Such a matrix $K$ of size $N\times N$
is said to satisfy the {\bf complementary low-rank property} if for
any level $\ell$, any node $A$ in $T_X$ at level $\ell$, and any node
$B$ in $T_\Omega$ at level $L-\ell$, the submatrix $K_{A,B}$, obtained
by restricting $K$ to the rows indexed by the points in $A$ and the
columns indexed by the points in $B$, is numerically low-rank, i.e.,
for a given precision $\epsilon$ there exists a low-rank approximation
of $K_{A,B}$ with the $2$-norm error bounded by $\epsilon$ and the
rank bounded polynomially in $\log N$ and $\log(1/\epsilon)$.  In many
applications, one can even show that the rank is only bounded
polynomially in $\log(1/\epsilon)$ and is independent of $N$. While it
is straightforward to generalize the concept of the complementary
low-rank property to a matrix with different row and column
dimensions, the following discussion is restricted to the square
matrices for simplicity.

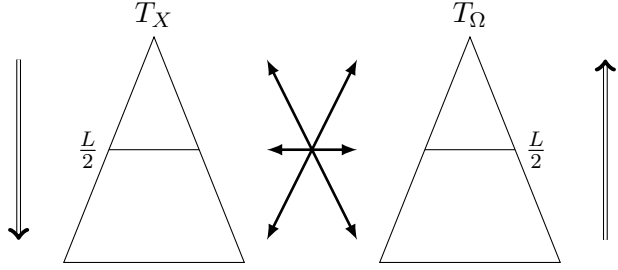
\begin{figure}
\centering
\begin{tikzpicture}[scale=0.3]

\draw [->,double,double distance=1pt] (-6,9) -> (-6,1);

\draw (0,10) -- (-4,0) -- (4,0) -- cycle;

\draw [->,double,double distance=1pt] (20,1) -> (20,9);

\draw (14,10) -- (10,0) -- (18,0) -- cycle;

\draw [latex-latex,line width=1pt] (9,5) -> (5,5);
\draw [latex-latex,line width=1pt] (9,9) -> (5,1);
\draw [latex-latex,line width=1pt] (9,1) -> (5,9);

\coordinate [label=above:$T_X$] (X) at (0,10);
\coordinate [label=above:$T_\Omega$] (X) at (14,10);

\draw (-2,5) -- (2,5);
\draw (12,5) -- (16,5);

\coordinate [label=right:$\frac{L}{2}$] (X) at (16,5);
\coordinate [label=left:$\frac{L}{2}$] (X) at (-2,5);
\end{tikzpicture}
\caption{Trees of the row and column indices.
    Left: $T_X$ for the row indices $X$.
    Right: $T_\Omega$ for the column indices $\Omega$.
    The interaction between $A\in T_X$ and $B\in T_\Omega$
    starts at the root of $T_X$ and the leaves of $T_\Omega$. }
\label{fig:domain-tree-BF}
\end{figure}

A simple yet important example is the Fourier matrix $K$ of size
$N\times N$, where
\begin{align*}
  & X = \Omega = \{0,\ldots, N-1\},\\
  & K = \left(\exp(2\pi \imath jk/N) \right)_{0\le j,k< N}.
\end{align*}
Here the trees $T_{X}$ and $T_{\Omega}$ are generated by bisecting the
sets $X$ and $\Omega$ recursively. Both trees have the same depth
$L=\log_2 N$. For each pair of nodes $A \in T_{X}$ and $B\in
T_{\Omega}$ with $A$ at level $\ell$ and $B$ at level $L-\ell$, the
numerical rank of the submatrix $K_{A,B}$ for a fixed precision
$\epsilon$ is bounded by a number that is independent of $N$ and
scales linearly with respect to $\log(1/\epsilon)$ \cite{Butterfly2}.

\begin{figure}[ht!]
  \begin{center}
    \begin{tabular}{ccccc}
      \includegraphics[height=1.1in]{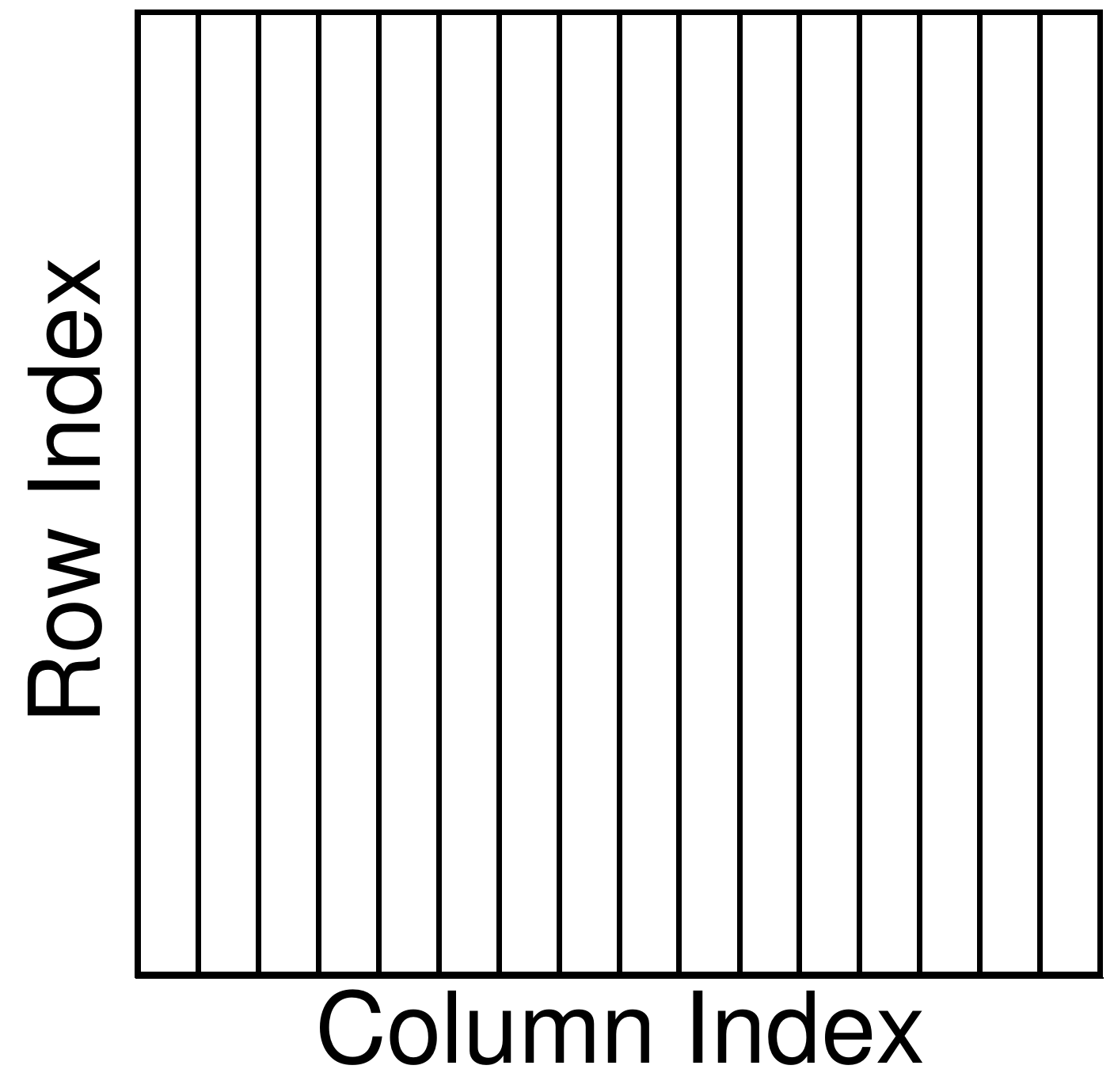}&
      \includegraphics[height=1.1in]{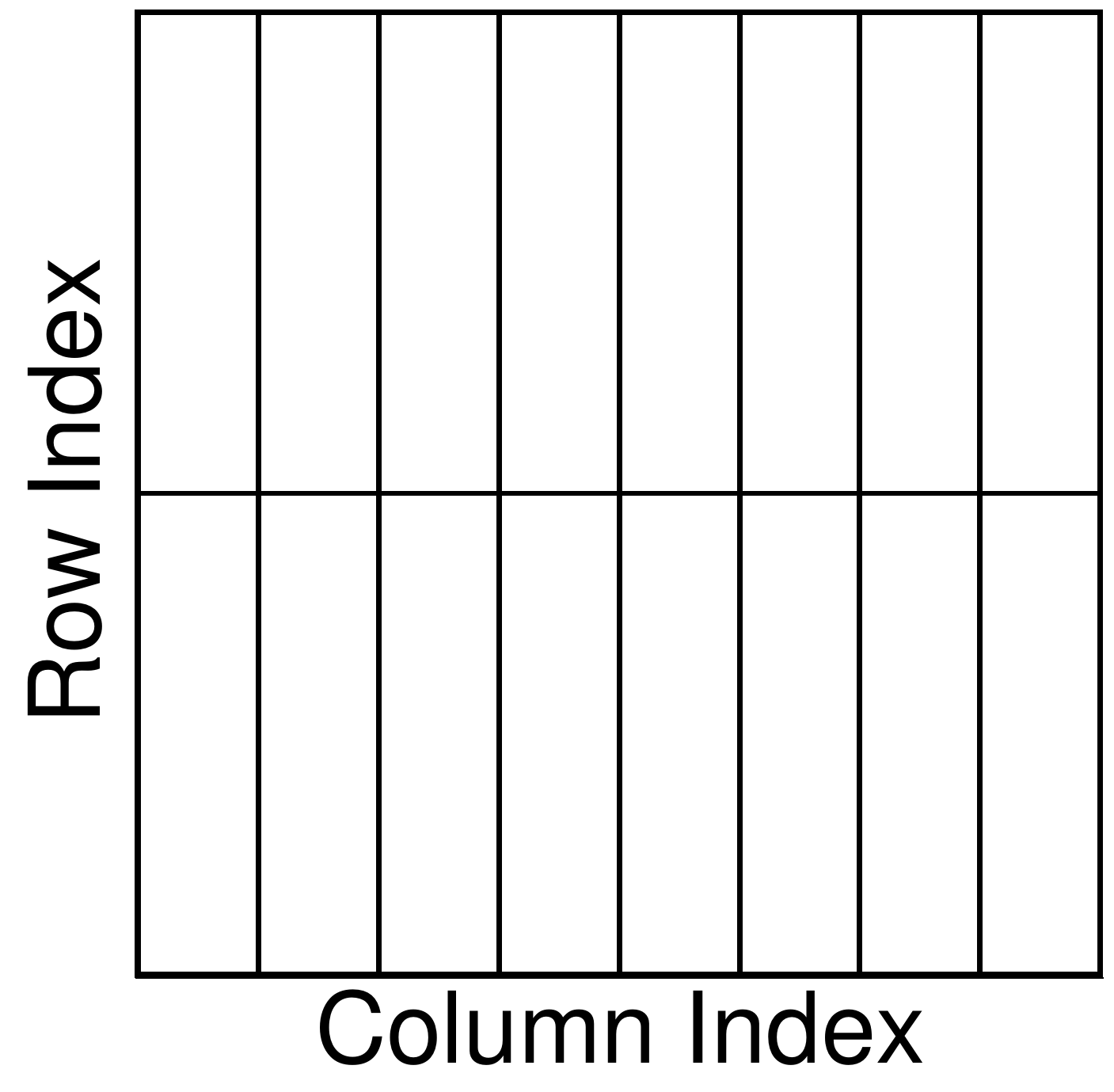}&
      \includegraphics[height=1.1in]{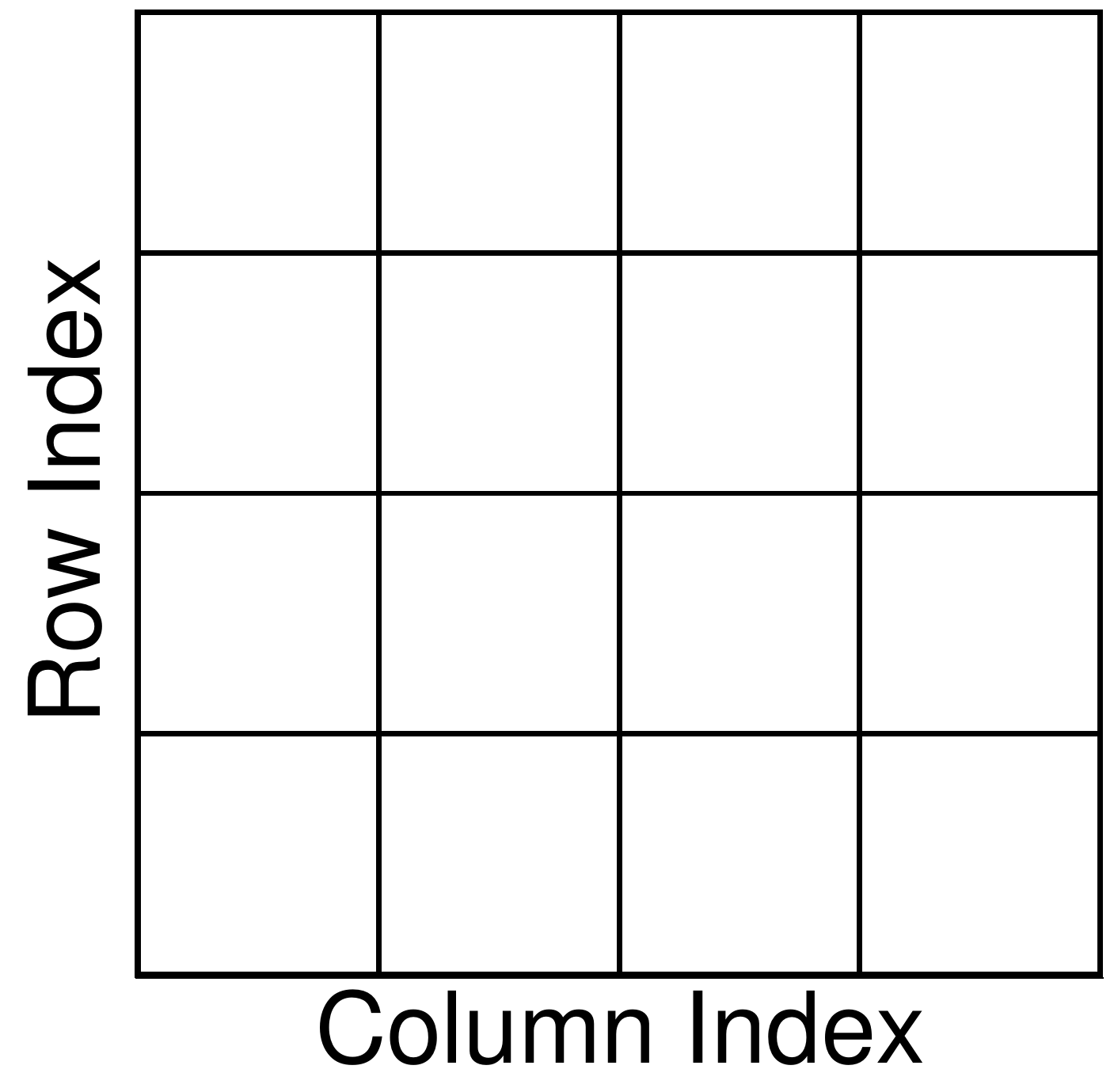}&
      \includegraphics[height=1.1in]{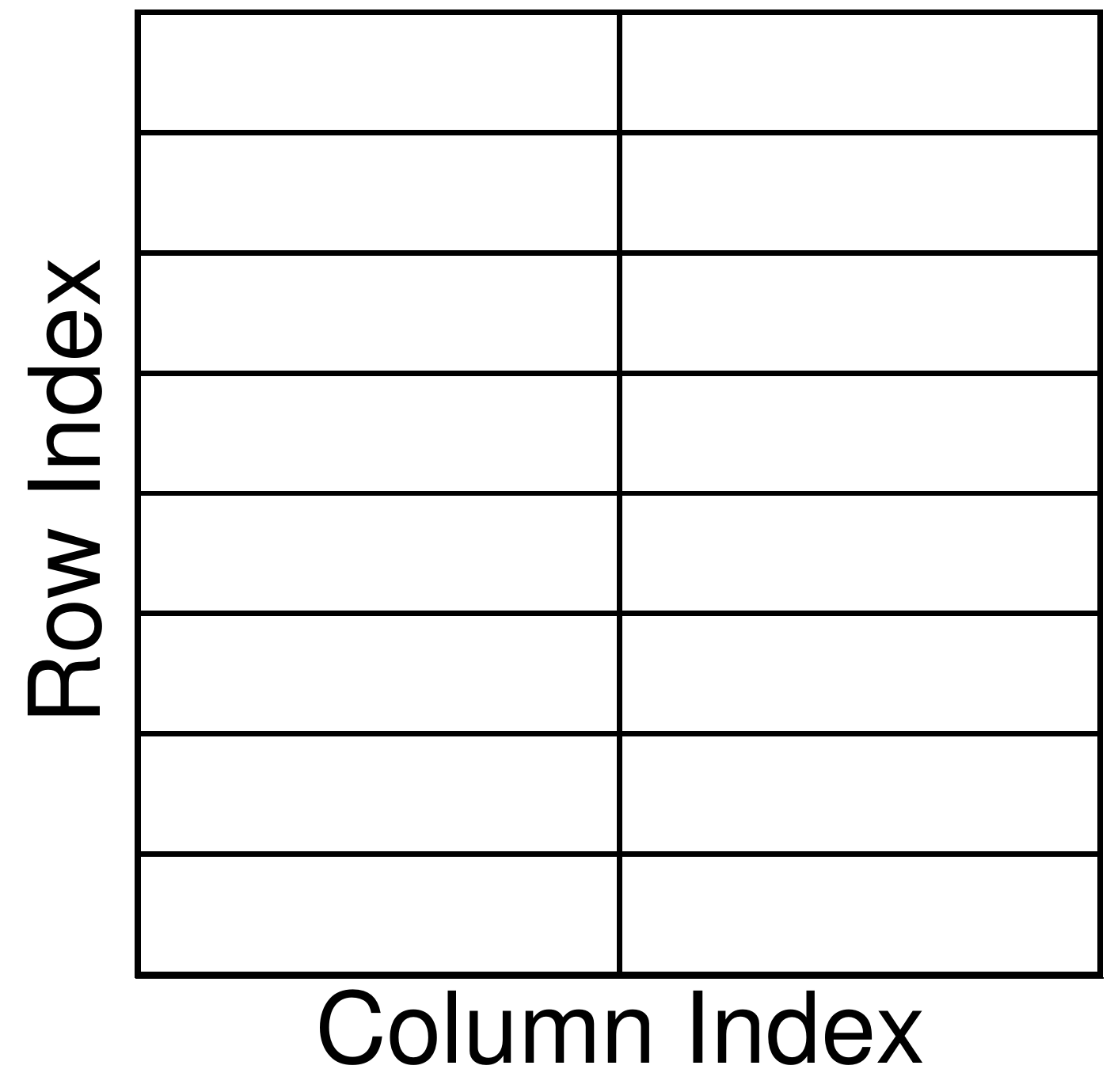}&
      \includegraphics[height=1.1in]{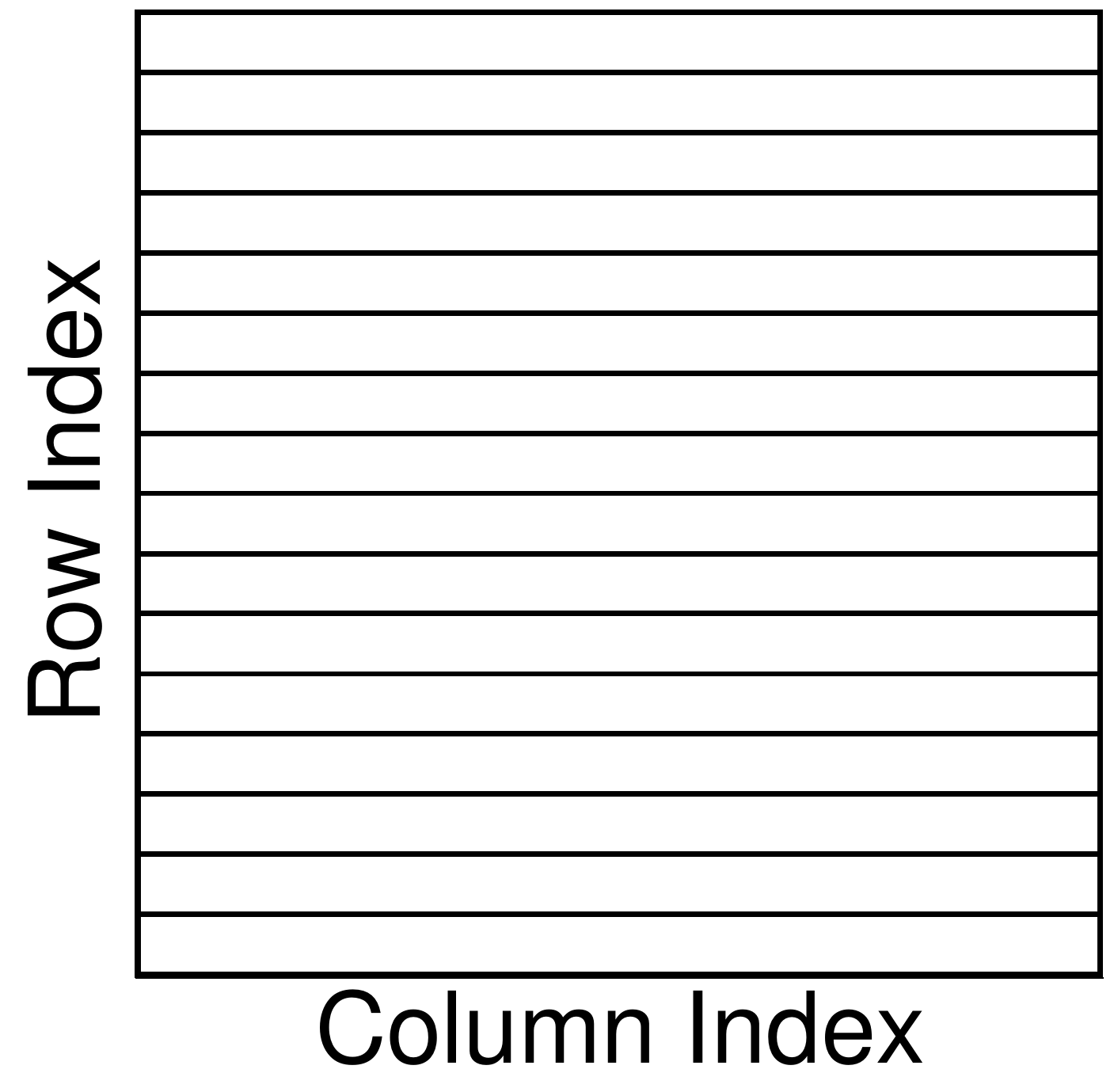}
    \end{tabular}
  \end{center}
  \caption{Hierarchical decomposition of the row and column indices
    of a $16\times 16$ matrix.  The trees $T_{X}$ and $T_{\Omega}$
    have roots containing $16$ column and row indices and leaves
    containing a single column and row index.  The rectangles above
    indicate the submatrices satisfying the complementary low-rank property.}
\label{fig:submatrices}
\end{figure}

For complementary low-rank matrices, the matrix-vector multiplication
can be carried out efficiently via the {\bf butterfly algorithm},
which was initially proposed in \cite{Butterfly1} and later extended
in \cite{Butterfly2}. For a general matrix $K$ of this type, the
butterfly algorithm consists of two stages: the off-line stage and the
on-line stage. In the off-line stage, it conducts simultaneously a top
down traversal of $T_{X}$ and a bottom up traversal of $T_{\Omega}$
(see Figure \ref{fig:domain-tree-BF} for an interpretation of data
flows) to recursively compress all complementary low-rank submatrices
(see Figure \ref{fig:submatrices} for an example of necessary
submatrices).  This typically takes $O(N^2)$ operations
\cite{Butterfly2,wavemoth} for a general complementary low-rank matrix $K$. In
the on-line stage, the butterfly algorithm then evaluates $u=Kg$ for a
given input vector $g\in\bbC^N$ in $O(N\log N)$ operations. While the
on-line application cost is essentially linear, the $O(N^2)$ off-line
precomputation cost appears to be a major bottleneck for many
calculations. For certain complementary low-rank matrices, such as the
ones obtained from the Fourier integral operators ({FIOs})
\cite{FIO09,FIO14,FIO13}, the sparse Fourier transforms
\cite{Butterfly5}, and the numerical solutions of acoustic wave
equations \cite{FIO12}, the off-line precomputation cost can be
reduced to nearly linear or even totally eliminated. However, in all
these cases, the reduction heavily relies on strong assumptions on the
analytic properties of the kernel function of $K$. When such detailed
information is not available, we are then forced to fall back on the
$O(N^2)$ off-line precomputation algorithm.

\subsection{Motivations and significance}

A natural question is whether it is still possible to reduce the cost
of the precomputation stage if the analytic properties of the kernel
are not accessible. The following two cases are quite common in
applications:
\begin{enumerate}
\item Only black-box routines for computing $Kg$ and $K^*g$ in
  $O(N\log N)$ operations are given.
\item Only a black-box routine for evaluating any entry of the matrix
  $K$ in $O(1)$ operations is given.
\end{enumerate}

To answer this question, this paper proposes the {\bf butterfly
  factorization}, which represents $K$ as a product of $L+3$ sparse
matrices:
\begin{equation}
\label{eqn:Kfactors}
K \approx U^{L}G^{L-1}\cdots G^{h} M^h (H^{h})^* \cdots (H^{L-1})^* (V^L)^*,
\end{equation}
where the depth $L=\O(\log N)$ of $T_X$ and $T_\Omega$ is assumed to be even,
$h=L/2$ is a middle level index, and all factors are sparse matrices
with $\O(N)$ nonzero entries.

The construction of the butterfly factorization proceeds as follows in
two stages. The first stage is to construction a preliminary middle
level factorization that is associated with the middle level of $T_X$
and $T_\Omega$
\begin{equation}
  \label{eqn:intFac}
  K\approx U^hM^h (V^h)^*,
\end{equation}
where $U^h$ and $V^h$ are block diagonal matrices and $M^h$ is a
weighted permutation matrix. In the first case, this is achieved by
applying $K$ to a set of $O(N^{1/2})$ structured random vectors and
then applying the randomized singular value decomposition (SVD) to the
result. This typically costs $O(N^{3/2}\log N)$ operations. In the
second case, \eqref{eqn:intFac} is built via the randomized sampling
method proposed in \cite{Butterfly4,symbol3} for computing approximate
SVDs. This randomized sampling needs to make the assumption that the
columns and rows of middle level blocks of $K$ to be incoherent with
respect to the delta functions and it typcially takes only
$\O(N^{3/2})$ operations in practice.

Once the middle level factorization \eqref{eqn:intFac} is available,
the second stage is a sequence of truncated SVDs that further
factorize each of $U^h$ and $V^h$ into a sequence of sparse matrices,
resulting in the final factorization \eqref{eqn:Kfactors}. The
operation count of this stage is $\O(N^{3/2})$ and the total memory
complexity for constructing butterfly factorization is $\O(N^{3/2})$.

When the butterfly factorization \eqref{eqn:Kfactors} is constructed,
the cost of applying $K$ to a given vector $g\in\bbC^N$ is $\O(N\log
N)$ because \eqref{eqn:Kfactors} is a sequence of $O(\log N)$ sparse
matrices, each with $O(N)$ non-zero entries.  Although we shall limit
our discussion to one-dimensional problems in this paper, the proposed
butterfly factorization, along with its construction algorithm, can be
easily generalized to higher dimensions.

This work is motivated by problems that require repeated
applications of a butterfly algorithm.
In several applications, such as inverse
scattering \cite{InvRadon,InvFIO} and fast spherical harmonic
transform ({SHT}) \cite{Special2}, the butterfly algorithm is called
repeatedly either in an iterative process of minimizing some
regularized objective function or to a large set of different input
vectors. Therefore, it becomes important to reduce the
constant prefactor of the butterfly algorithm to save actual runtime.
For example in \cite{FIO09},
{Chebyshev} interpolation is applied to recover low-rank structures
of submatrices with a sufficiently large number of interpolation points.
The recovered rank is far from the optimum.
Hence, the prefactor of the corresponding butterfly algorithm
in \cite{FIO09} is large.
The butterfly factorization can further compress this butterfly algorithm to obtain
nearly optimal low-rank approximations resulting in a much smaller prefactor,
as will be shown in the numerical results. Therefore, it is
more efficient to construct the butterfly factorization using this
butterfly algorithm and then apply the butterfly factorization
repeatedly.
In this sense, the butterfly factorization can be viewed as
a compression of certain butterfly algorithms.

Another important application is
the computation of a composition of several {FIOs}.
A direct method to construct the composition takes $\O(N^3)$ operations,
while the butterfly factorization provides a data-sparse representation
of this composition in $\O(N^{3/2}\log N)$ operations,
once the fast algorithm for applying each {FIO} is available.
After the construction, the application of the butterfly factorization
is independent of the number of {FIOs} in the composition,
which is significant when the number of {FIOs} is large.

Recently, there has also been a sequence of papers on recovering a
structured matrix via applying it to (structured) random vectors. For
example, the randomized SVD algorithms \cite{Rec1,Rec2,Rec3} recover a
low-rank approximation to an unknown matrix when it is numerically
low-rank. The work in \cite{HSSMatrix} constructs a sparse
representation for an unknown {HSS}
matrix. More recently, \cite{HMatrix} considers the more general
problem of constructing a sparse representation of an unknown
$\mathcal{H}$-matrix. To our best knowledge, the present work is the
first to address such matrix recovery problem if the unknown matrix
satisfies the complementary low-rank property.

\subsection{Content}

The rest of this paper is organized as follows. Section \ref{sec:prep}
briefly reviews some basic tools that shall be used repeatedly in
Sections \ref{sec:butterflyfac}. Section \ref{sec:butterflyfac}
describes in detail the butterfly factorization and its construction
algorithm. In Section \ref{sec:results}, numerical examples are
provided to demonstrate the efficiency of the proposed algorithms.
Finally, Section \ref{sec:conclusion} lists several directions for
future work.

\section{Preliminaries}
\label{sec:prep}

For a matrix $Z \in \bbC^{m\times n}$, we define a rank-$r$
approximate singular value decomposition (SVD) of $Z$ as
\[
Z \approx U_0 \Sigma_0 V_0^*,
\]
where $U_0\in \bbC^{m\times r}$ is unitary, $\Sigma_0\in \bbR^{r\times
  r}$ is diagonal, and $V_0\in \bbC^{n\times r}$ is unitary. A
straightforward method to obtain the optimal rank-$r$ approximation of
$Z$ is to compute its truncated SVD, where $U_0$ is the matrix with
the first $r$ left singular vectors, $\Sigma_0$ is a diagonal matrix
with the first $r$ singular values in decreasing order, and $V_0$ is
the matrix with the first $r$ right singular vectors.

A typical computation of the truncated {SVD} of $Z$ takes
$\O(mn\min(m,n))$ operations, which can be quite expensive when $m$
and $n$ are large. Therefore, a lot of research has been devoted to
faster algorithms for computing approximate SVDs, especially for
matrices with fast decaying singular values. In Sections
\ref{sec:RandSVD} and \ref{sec:RandSampling}, we will introduce two
randomized algorithms for computing approximate SVDs for numerically
low-rank matrices $Z$: the first one \cite{Rec1} is based on applying
the matrix to random vectors while the second one
\cite{Butterfly4,symbol3} relies on sampling the matrix entries
randomly.

Once an approximate SVD $Z \approx U_0 \Sigma_0 V_0^*$ is computed, it
can be written in several equivalent ways, each of which is convenient
for certain purposes. First, one can write
\[
Z \approx U S V^*,
\]
where
\begin{equation}
  \label{eq:lowrankSVD}
  U=U_0\Sigma_0,\, S = \Sigma_0^{-1} \text{ and } V^*=\Sigma_0V_0^*.
\end{equation}
This construction is analogous to the well-known {CUR}
decomposition \cite{mahoney2009cur} in the sense that the left and
right factors in both factorization methods inherit similar singular
values of the original numerical low-rank matrix. Here, the middle
matrix $S$ in \eqref{eq:lowrankSVD} can be carefully constructed to
ensure numerical stability, since the singular values in $\Sigma_0$
can be computed to nearly full relative precision.

As we shall see, sometimes it is also convenient to write the
approximation as
\[
Z\approx UV^*
\]
where
\begin{equation}
  \label{eq:lowrankSVD2}
  U=U_0 \text{ and } V^*=\Sigma_0V_0^*,
\end{equation}
or
\begin{equation}
  \label{eq:lowrankSVD3}
  U=U_0 \Sigma_0\text{ and }
  V^*=V_0^*.
\end{equation}
Here, one of the factors $U$ and $V$ share the singular values of $Z$.

\subsection{SVD via random matrix-vector multiplication}
\label{sec:RandSVD}

One popular approach is the randomized algorithm in \cite{Rec1} that
reduces the cubic complexity to $\O(rmn)$ complexity. We briefly
review this following \cite{Rec1} for constructing a rank-$r$
approximation SVD $Z\approx U_0 \Sigma_0 V_0^*$ below.
\begin{algo}{Randomized {SVD}}
\begin{enumerate}
\item
  Generate two tall skinny random Gaussian matrices $R_{col} \in
  \bbC^{n\times (r+p)}$ and $R_{row}\in \bbC^{m\times
    (r+p)}$, where $p=\O(1)$ is an additive oversampling parameter that increases
  the approximation accuracy.
\item
  Apply the pivoted QR factorization to $ZR_{col}$ and let $Q_{col}$
  be the matrix of the first $r$ columns of the $Q$ matrix. Similarly,
  apply the pivoted QR factorization to $Z^*R_{row}$ and let $Q_{row}$
  be the matrix of the first $r$ columns of the $Q$ matrix.
\item
  Generate a tiny middle matrix $M=Q_{col}^*ZQ_{row}$ and compute its
  rank-$r$ truncated {SVD}: $M \approx U_M \Sigma_M V_M^*$.
\item
  Let $U_0 = Q_{col} U_M$, $\Sigma_0 = \Sigma_M$, and $V_0^* = V_M^*
  Q_{row}^*$.  Then $Z \approx U_0\Sigma_0V_0^*$.
\end{enumerate}
\end{algo}
The dominant complexity comes from the application of $Z$ to $\O(r)$
random vectors. If fast algorithms for applying $Z$ are available, the
quadratic complexity can be further reduced.

Once the approximate {SVD} of $Z$ is ready, the equivalent forms in
\eqref{eq:lowrankSVD}, \eqref{eq:lowrankSVD2}, and
\eqref{eq:lowrankSVD3} can be constructed easily. Under the condition
that the singular values of $Z$ decay sufficiently rapidly, the
approximation error of the resulting rank-$r$ is nearly optimal with
an overwhelming probability. Typically, the additive over-sampling
parameter $p=5$ is sufficient to obtain an accurate rank-$r$
approximation of $Z$.


For most applications, the goal is to construct a low-rank
approximation up to a fixed relative precision $\epsilon$, rather than
a fixed rank $r$. The above procedure can then be embedded into an
iterative process that starts with a relatively small $r$, computes a
rank-$r$ approximation, estimates the error probabilistically, and
repeats the steps with doubled rank $2r$ if the error is above the
threshold $\epsilon$ \cite{Rec1}.

\subsection{SVD via random sampling}
\label{sec:RandSampling}
The above algorithm relies only on the product of the matrix $Z\in
\bbC^{m\times n}$ or its transpose with given random vectors.  If one
is allowed to access the individual entries of $Z$, the following
randomized sampling method for low-rank approximations introduced in
\cite{Butterfly4,symbol3} can be more efficient. This method only
visits $\O(r)$ columns and rows of $Z$ and hence only requires
$\O(r(m+n))$ operations and memory.

Here, we adopt the standard notation for a submatrix: given a row
index set $I$ and a column index set $J$, $Z_{I,J}=Z(I,J)$ is the
submatrix with entries from rows in $I$ and columns in $J$; we also
use $``:"$ to denote the entire columns or rows of the matrix, i.e.,
$Z_{I,:} = Z(I,:)$ and $Z_{:,J} = Z(:,J)$.  With these handy
notations, we briefly introduce the randomized sampling algorithm to
construct a rank-$r$ approximation of $Z\approx U_0\Sigma_0V_0^*$.

\begin{algo}{Randomized sampling for low-rank approximation}
\begin{enumerate}
\item
  Let $\Pi_{col}$ and $\Pi_{row}$ denote the important columns and rows of $Z$
  that are used to form the column and row bases.
  Initially $\Pi_{col} = \emptyset$ and $\Pi_{row} = \emptyset$.
\item
  Randomly sample $rq$ rows and denote their indices by $S_{row}$.  Let $I
  = S_{row}\cup \Pi_{row}$. Here $q=\O(1)$ is a multiplicative oversampling
  parameter. Perform a pivoted {QR} decomposition of $Z_{I,:}$ to get
  \[
  Z_{I,:} P = QR,
  \]
  where $P$ is the resulting permutation matrix and $R=(r_{ij})$ is an
  $\O(r)\times n$ upper triangular matrix. Define the important column
  index set $\Pi_{col}$ to be the first $r$ columns picked within the
  pivoted QR decomposition.
\item
  Randomly sample $rq$ columns and denote their indices by $S_{col}$.
  Let $J = S_{col}\cup \Pi_{col}$. Perform a pivoted {LQ} decomposition of
  $Z_{:,J}$ to get
  \[
  P Z_{:,J} = LQ,
  \]
  where $P$ is the resulting permutation matrix and $L=(l_{ij})$ is an
  $m\times \O(r)$ lower triangular matrix.  Define the important
  row index set $\Pi_{row}$ to be the first $r$ rows picked within the
  pivoted LQ decomposition.
\item
  Repeat steps 2 and 3 a few times to ensure $\Pi_{col}$ and $\Pi_{row}$
  sufficiently sample the important columns and rows of $Z$.
\item
  Apply the pivoted QR factorization to $Z_{:,\Pi_{col}}$ and let
  $Q_{col}$ be the matrix of the first $r$ columns of the $Q$ matrix.
  Similarly, apply the pivoted QR factorization to $Z_{\Pi_{row},:}^*$
  and let $Q_{row}$ be the matrix of the first $r$ columns of the $Q$
  matrix.
\item
  We seek a middle matrix $M$ such that $Z \approx Q_{col} M
  Q_{row}^*$. To solve this problem efficiently, we approximately
  reduce it to a least-squares problem of a smaller size.  Let
  $S_{col}$ and $S_{row}$ be the index sets of a few extra randomly
  sampled columns and rows.  Let $J = \Pi_{col} \cup S_{col}$ and $I =
  \Pi_{row} \cup S_{row}$.  A simple least-squares solution to the
  problem
  \[
  \min_{M} \norm{Z_{I,J} - (Q_{col})_{I,:} M (Q_{row}^*)_{:,J}}
  \]
  gives $M = (Q_{col})_{I,:}^\dagger Z_{I,J}
  (Q_{row}^*)_{:,J}^\dagger$, where $(\cdot)^\dagger$ stands for the
  pseudo-inverse.
\item
  Compute an {SVD} $M \approx U_M\Sigma_MV_M^*$. Then the low-rank
  approximation of $Z\approx U_0 S_0 V_0^*$ is given by
  \begin{equation}
    U_0 = Q_{col} U_M,\quad \Sigma_0 = \Sigma_M, \quad V_0^* = V_M^*Q_{row}^* .
  \end{equation}
\end{enumerate}
\end{algo}

We have not been able to quantify the error and success probability
rigorously for this procedure at this point.
On the other hand,
when the columns and rows of $K$ are incoherent with respect to ``delta
functions'' (i.e., vectors that have only one significantly larger entry),
this procedure works well in our numerical experiments.
Here, a vector $u$ is said to be incoherent with respect to a vector $v$ if
$\mu = |u^Tv|/(\norm{u}_2\norm{v}_2)$ is small.
In the typical implementation, the multiplicative oversampling parameter
$q$ is equal to $3$ and Steps 2 and 3 are iterated no more than three
times.
These parameters are empirically sufficient to achieve accurate
low-rank approximations and are used through out numerical examples
in Section \ref{sec:results}.

As we mentioned above, for most applications the goal is to construct a
low-rank approximation up to a fixed relative error $\epsilon$, rather
than a fixed rank. This process can also be embedded into an
iterative process to achieve the desired accuracy.

\section{Butterfly factorization}
\label{sec:butterflyfac}

This section presents the butterfly factorization algorithm for a
matrix $K\in\bbC^{N\times N}$. For simplicity let
$X=\Omega=\{1,\ldots, N\}$. The trees $T_X$ and $T_\Omega$ are
complete binary trees with $L=\log_2 N - O(1)$ levels. We assume that
$L$ is an even integer and the number of points in each leaf node of
$T_X$ and $T_\Omega$ is bounded by a uniform constant.

At each level $\ell$, $\ell=0,\dots,L$, we denote the $i$th node at
level $\ell$ in $T_X$ as $A^\ell_i$ for $i=0,1,\dots,2^\ell-1$ and the
$j$th node at level $L-\ell$ in $T_\Omega$ as $B^{L-\ell}_j$ for
$j=0,1,\dots,2^{L-\ell}-1$. These nodes naturally partition $K$ into
$\O(N)$ submatrices $K_{A^\ell_i,B^{L-\ell}_j}$. For simplicity, we
write $K^\ell_{i,j}:=K_{A^\ell_i,B^{L-\ell}_j}$, where the superscript
is used to indicate the level (in $T_X$). The butterfly factorization
utilizes rank-$r$ approximations of all submatrices $K^\ell_{i,j}$
with $r=\O(1)$.

The butterfly factorization of $K$ is built in two stages.  In the
first stage, we compute a rank-$r$ approximations of each submatrix
$K^h_{i,j}$ at the level $\ell=h=L/2$ and then organize them into an
initial factorization:
\[
K\approx U^h M^h (V^h)^*,
\]
where $U^h$ and $V^h$ are block diagonal matrices and $M^h$ is a
weighted permutation matrix. This is referred as the {\bf middle level
  factorization} and is described in detail in Section \ref
{sec:MiddleConstruction}.

In the second stage, we recursively factorize $U^\ell \approx
U^{\ell+1}G^{\ell}$ and $(V^\ell)^* \approx (H^{\ell})^*
(V^{\ell+1})^*$ for $\ell=h,h+1,\dots,L-1$, since $U^\ell$ and
$(V^{\ell})^*$ inherit the complementary low-rank property from $K$,
i.e., the low-rank property of $U^\ell$ comes from the low-rank
property of $K^\ell_{i,j}$ and the low-rank property of $V^\ell$
results from the one of $K^{L-\ell}_{i,j}$. After this recursive
factorization, one reaches at the the butterfly factorization of $K$
\begin{equation}
  K \approx U^{L}G^{L-1}\cdots G^{h} M^h (H^{h})^* \cdots (H^{L-1})^*
  (V^L)^*,
\end{equation}
where all factors are sparse matrices with $\O(N)$ nonzero entries.
We refer to this stage as the {\bf recursive factorization} and it is
discussed in detail in Section \ref{sec:RC}.


\subsection{Middle level factorization}
\label{sec:MiddleConstruction}

The first step of the middle level factorization is to compute a
rank-$r$ approximation to every $K^h_{i,j}$. Recall that we
consider one of the following two cases.
\begin{enumerate}
\item Only black-box routines for computing $Kg$ and $K^*g$ in
  $O(N\log N)$ operations are given.
\item Only a black-box routine for evaluating any entry of the matrix
  $K$ in $O(1)$ operations is given.
\end{enumerate}
The actual computation of this step proceeds differently depending on
which case is under consideration. Through the discussion, $m=2^h =
O(N^{1/2})$ is the number of nodes in the middle level $h=L/2$ and we
assume without loss of generality that $N/m$ is an integer.
\begin{itemize}
\item In the first case, the rank-$r$ approximation of each
  $K^h_{i,j}$ is constructed with the SVD algorithm via random
  matrix-vector multiplication in Section \ref{sec:RandSVD}. This
  requires us to apply $K^h_{i,j}$ and its adjoint to random Gaussian
  matrices of size $(N/m)\times(r+p)$, where $r$ is the desired rank
  and $p$ is an oversampling parameter. In order to take advantage of
  the fast algorithm for multiplying $K$, we construct a matrix $C$ of
  size $N \times m(r+p)$. $C$ is partitioned into an $m\times m$
  blocks with each block $C_{ij}$ for $i,j=0,1,\dots,m-1$ of size
  $(N/m)\times (r+p)$. In additional, $C$ is block-diagonal and its
  diagonal blocks are random Gaussian matrices.  This is equivalent to
  applying each $K^h_{i,j}$ to the same random Gaussian matrix
  $C_{jj}$ for all $i$. We then use the fast algorithm to apply $K$ to
  each column of $C$ and store the results. Similarly, we form another
  random block diagonal matrix $R$ similar to $C$ and use the fast
  algorithm of applying $K^*$ to $R$. This is equivalent to applying
  each $(K^h_{i,j})^*$ to an $(N/m)\times (r+p)$ Gaussian random
  matrix $R_{ii}$ for all $j=0,1,\ldots, m-1$. With $K^h_{i,j}C_{jj}$ and
  $(K^h_{i,j})^*R_{ii}$ ready, we can compute the rank-$r$ approximate
  SVD of $K^h_{i,j}$ following the procedure described in Section
  \ref{sec:RandSVD}.
\item In the second case, it is assumed that an arbitrary entry of $K$
  can be calculated in $\O(1)$ operations. We simply apply the SVD
  algorithm via random sampling in Section \ref{sec:RandSampling} to
  each $K^h_{i,j}$ to construct a rank-$r$ approximate SVD.
\end{itemize}
In either case, once the approximate SVD of $K^h_{i,j}$ is ready, it
is transformed in the form
\[
K^h_{i,j} \approx U^h_{i,j}S^h_{i,j}(V^h_{j,i})^*
\]
following \eqref{eq:lowrankSVD}. We would like to emphasize that the
columns of $U^h_{i,j}$ and $V^h_{j,i}$ are scaled with the singular
values of the approximate SVD so that they keep track of the
importance of these columns in approximating $K^h_{i,j}$.

\begin{figure}[htp]
\begin{minipage}{\textwidth}
\centering
\resizebox{3.5cm}{!}{
\begin{tikzpicture}[baseline=-0.5ex]
      \tikzset{every left delimiter/.style={xshift=-1ex},every right delimiter/.style={xshift=1ex}}
      \matrix (mat) [matrix of math nodes, left delimiter=(, right delimiter=)] {
\draw[fill=gray] (0,0) rectangle (32,32);
\\
      };
\end{tikzpicture}
}
$\approx$
\resizebox{3.5cm}{!}{
\begin{tikzpicture}[baseline=-0.5ex]
      \tikzset{every left delimiter/.style={xshift=-1ex},every right delimiter/.style={xshift=1ex}}
      \matrix (mat) [matrix of math nodes, left delimiter=(, right delimiter=)] {
      \draw;
\foreach \i in {0,0.5,...,3.5} {
    \foreach \j in {0,4,...,28} {
        \draw[fill=gray] (32-\i-\j,\j) rectangle (32-\i-\j-0.5,\j+4);
    }
}
\\
      };
\end{tikzpicture}
}
\resizebox{3.5cm}{!}{
\begin{tikzpicture}[baseline=-0.5ex]
      \tikzset{every left delimiter/.style={xshift=-1ex},every right delimiter/.style={xshift=1ex}}
      \matrix (mat) [matrix of math nodes, left delimiter=(, right delimiter=)] {
      \draw;
\foreach \i in {0,0.5,...,3.5}{
    \foreach \j in {0,0.5,...,3.5}{
        \draw[fill=gray] (8*\i+\j,16-8*\j-\i) rectangle (8*\i+\j+0.5,15.5-8*\j-\i);
    }
}
\\
      };
\end{tikzpicture}
}
\resizebox{3.5cm}{!}{
\begin{tikzpicture}[baseline=-0.5ex]
      \tikzset{every left delimiter/.style={xshift=-1ex},every right delimiter/.style={xshift=1ex}}
      \matrix (mat) [matrix of math nodes, left delimiter=(, right delimiter=)] {
            \draw;
\foreach \i in {0,0.5,...,3.5} {
    \foreach \j in {0,4,...,28} {
        \draw[fill=gray] (\j,32-\i-\j) rectangle (\j+4,32-\i-\j-0.5);
    }
}
\\
      };
\end{tikzpicture}
}\\
\end{minipage}
\caption{The middle level factorization of
    a $64\times 64$ complementary low-rank matrix $K\approx U^3M^3(V^3)^*$
    assuming $r=1$.
    Grey blocks indicate nonzero blocks.
    $U^3$ and $V^3$ are block-diagonal matrices with $8$ blocks.
    The diagonal blocks of $U^3$ and $V^3$ are assembled
    according to Equation \eqref{eq:expression-UV}
    as indicated by black rectangles.
    $M^3$ is a $8\times 8$ block matrix with each block
    $M^3_{i,j}$ itself an $8\times 8$ block matrix
    containing diagonal weights matrix on the $(j,i)$ block.}
\label{fig:compression-1D}
\end{figure}
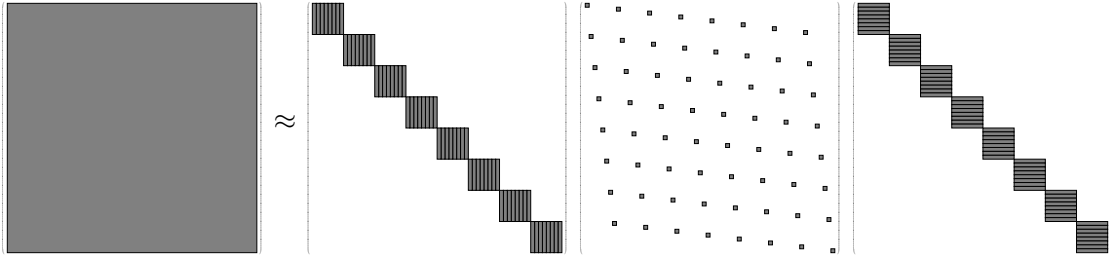

After calculating the approximate rank-$r$ factorization of each
$K^h_{i,j}$, we assemble these factors into three block matrices
$U^h$, $M^h$ and $V^h$ as follows:
\begin{equation}\label{eq:KeqUMV}
  \begin{split}
    K \approx &
    \begin{pmatrix}
      U^h_{0,0}S^h_{0,0}(V^h_{0,0})^* & U^h_{0,1}S^h_{0,1}(V^h_{1,0})^* &
      \cdots & U^h_{0,m-1}S^h_{0,m-1}(V^h_{m-1,0})^*\\
      U^h_{1,0}S^h_{1,0}(V^h_{0,1})^* & U^h_{1,1}S^h_{1,1}(V^h_{1,1})^* &
              & U^h_{1,m-1}S^h_{1,m-1}(V^h_{m-1,1})^*\\
      \vdots & & \ddots & \\
      U^h_{m-1,0}S^h_{m-1,0}(V^h_{0,m-1})^* & U^h_{m-1,1}S^h_{m-1,1}(V^h_{1,m-1})^* &
              & U^h_{m-1,m-1}S^h_{m-1,m-1}(V^h_{m-1,m-1})^*
    \end{pmatrix}\\
    = &
    \begin{pmatrix}
      U^h_{0} & & &\\
      & U^h_{1} & &\\
      & & \ddots &\\
      & & & U^h_{m-1}
    \end{pmatrix}
    \begin{pmatrix}
      M^h_{0,0} & M^h_{0,1} & \cdots & M^h_{0,m-1}\\
      M^h_{1,0} & M^h_{1,1} & & M^h_{1,m-1}\\
      \vdots & & \ddots & \\
      M^h_{m-1,0} & M^h_{m-1,1} & & M^h_{m-1,m-1}
    \end{pmatrix}
    \begin{pmatrix}
      (V_{0}^h)^* & & &\\
      & (V_{1}^h)^* & &\\
      & & \ddots &\\
      & & & (V_{m-1}^h)^*
    \end{pmatrix}\\
    = & U^h M^h (V^h)^*,\\
  \end{split}
\end{equation}
where
\begin{equation}
  U_{i}^h=
  \begin{pmatrix}
    U^h_{i,0} & U^h_{i,1} & \cdots & U^h_{i,m-1}
  \end{pmatrix}
  \in \bbC^{(N/m) \times m r},
  \quad
  V^h_{j}=
  \begin{pmatrix}
    V^h_{j,0} & V^h_{j,1} & \cdots & V^h_{j,m-1}
  \end{pmatrix}
  \in \bbC^{(N/m) \times m r},
  \label{eq:expression-UV}
\end{equation}
and $M^h \in \bbC^{(m^2 r) \times (m^2 r)}$ is a weighted permutation
matrix. Each submatrix $M^h_{i,j}$ is itself an $m\times m$ block
matrix with block size $r\times r$ where all blocks are zero except
that the $(j,i)$ block is equal to the diagonal matrix $S_{i,j}^h$.
It is obvious that there are only $\O(N)$ nonzero entries in $M^h$.
See Figure~\ref{fig:compression-1D} for an example of a middle level
factorization of a $64\times 64$ matrix with $r=1$.

\subsection{Recursive factorization}
\label{sec:RC}

In this section, we will recursively factorize
\begin{equation}
  \label{eqn:UF}
  U^\ell \approx U^{\ell+1}G^{\ell}
\end{equation}
for $\ell=h,h+1,\dots,L-1$ and
\begin{equation}
  \label{eqn:VF}
  (V^\ell)^* \approx (H^{\ell})^* (V^{\ell+1})^*
\end{equation}
for $\ell=h,h+1,\dots,L-1$.  After these recursive factorizations, we
can obtain the following butterfly factorization by substituting these
factorizations into \eqref{eq:KeqUMV}:
\begin{equation}
  K \approx U^{L}G^{L-1}\cdots G^{h} M^h (H^{h})^* \cdots (H^{L-1})^* (V^L)^*.
\end{equation}

\subsubsection{Recursive factorization of $U^h$}
\label{sec:RCU}

Each factorization at level $\ell$ in \eqref{eqn:UF} results from the
low-rank property of $K^\ell_{i,j}$ for $\ell\geq L/2$.  When
$\ell=h$, recall that
\begin{equation*}
U^h=
\begin{pmatrix}
U^h_{0} & & &\\
 & U^h_{1} & &\\
 & & \ddots &\\
 & & & U^h_{m-1}
\end{pmatrix}
\end{equation*}
and
\begin{equation*}
U_{i}^h=
\begin{pmatrix}
U^h_{i,0} & U^h_{i,1} & \cdots & U^h_{i,m-1}
\end{pmatrix}
\end{equation*}
with each $U^h_{i,j}\in \bbC^{(N/m)\times r}$.  We split $U^h_i$ and
each $U^h_{i,j}$ into halves by row, i.e.,
\begin{equation*}
U^h_i=
\begin{pmatrix}
U^{h,t}_{i}\\
\midrule
U^{h,b}_{i}
\end{pmatrix}
\text{  and   }
U^h_{i,j}=
\begin{pmatrix}
U^{h,t}_{i,j}\\
\midrule
U^{h,b}_{i,j}
\end{pmatrix},
\end{equation*}
where the superscript $t$ denotes the top half
and $b$ denotes the bottom half of a matrix.
Then we have
\begin{equation}
\label{eqn:half}
U^h_{i}=
\begin{pmatrix}
U^{h,t}_{i,0} &  U^{h,t}_{i,1} & \dots &U^{h,t}_{i,m-1}\\
\midrule
U^{h,b}_{i,0} & U^{h,b}_{i,1} &\dots &U^{h,b}_{i,m-1}
\end{pmatrix}.
\end{equation}
Notice that, for each $i=0,1,\dots, m-1$ and $j=0,1,\dots, m/2-1$, the columns of
\begin{equation}
  \begin{pmatrix}
    U^{h,t}_{i,2j} &  U^{h,t}_{i,2j+1}
  \end{pmatrix}
  \text{   and    }
  \begin{pmatrix}
    U^{h,b}_{i,2j} & U^{h,b}_{i,2j+1}
  \end{pmatrix}
  \label{eq:blocks}
\end{equation}
in \eqref{eqn:half} are in the column space of $K^{h+1}_{2i,j}$ and
$K^{h+1}_{2i+1,j}$, respectively.  By the complementary low-rank
property of the matrix $K$, $K^{h+1}_{2i,j}$ and $K^{h+1}_{2i+1,j}$
are numerical low-rank.  Hence $\left( U^{h,t}_{i,2j} U^{h,t}_{i,2j+1}
\right)$ and $\left( U^{h,b}_{i,2j} U^{h,b}_{i,2j+1} \right)$ are
numerically low-rank matrices in $\bbC^{(N/2m)\times 2r}$. Compute
their rank-$r$ approximations by the standard truncated {SVD},
transform it into the form of \eqref{eq:lowrankSVD3} and denote them
as
\begin{equation}
  \begin{pmatrix}
    U^{h,t}_{i,2j} &  U^{h,t}_{i,2j+1}
  \end{pmatrix}
  \approx U^{h+1}_{2i,j} G^{h}_{2i,j}
  \text{   and    }
  \begin{pmatrix}
    U^{h,b}_{i,2j} & U^{h,b}_{i,2j+1}
  \end{pmatrix}
  \approx U^{h+1}_{2i+1,j} G^{h}_{2i+1,j}
  \label{eqn:RSh}
\end{equation}
for $i=0,1,\dots,m-1$ and $j=0,1,\dots,m/2-1$.  The matrices in
\eqref{eqn:RSh} can be assembled into two new sparse matrices, such
that
\begin{equation*}
U^h \approx U^{h+1}G^{h}=
\begin{pmatrix}
U^{h+1}_{0} & & &\\
 & U^{h+1}_{1} & &\\
 & & \ddots &\\
 & & & U^{h+1}_{2m-1}
\end{pmatrix}
\begin{pmatrix}
G^{h}_{0} & & &\\
 & G^{h}_{1} & &\\
 & & \ddots &\\
 & & & G^{h}_{m-1}
\end{pmatrix},
\end{equation*}
where
\begin{equation*}
U_{i}^{h+1}=
\begin{pmatrix}
U^{h+1}_{i,0} & U^{h+1}_{i,1} & \cdots & U^{h+1}_{i,m/2-1}
\end{pmatrix}
\end{equation*}
for $i=0,1,\dots,2m-1$, and
\begin{equation*}
G^{h}_{i}=
\begin{pmatrix}
G^{h}_{2i,0} & & &\\
 & G^{h}_{2i,1} & &\\
 & & \ddots &\\
 & & & G^{h}_{2i,m/2-1}\\
\midrule
G^{h}_{2i+1,0} & & &\\
 & G^{h}_{2i+1,1} & &\\
 & & \ddots &\\
 & & & G^{h}_{2i+1,m/2-1}
\end{pmatrix}
\end{equation*}
for $i=0,1,\dots,m-1$.

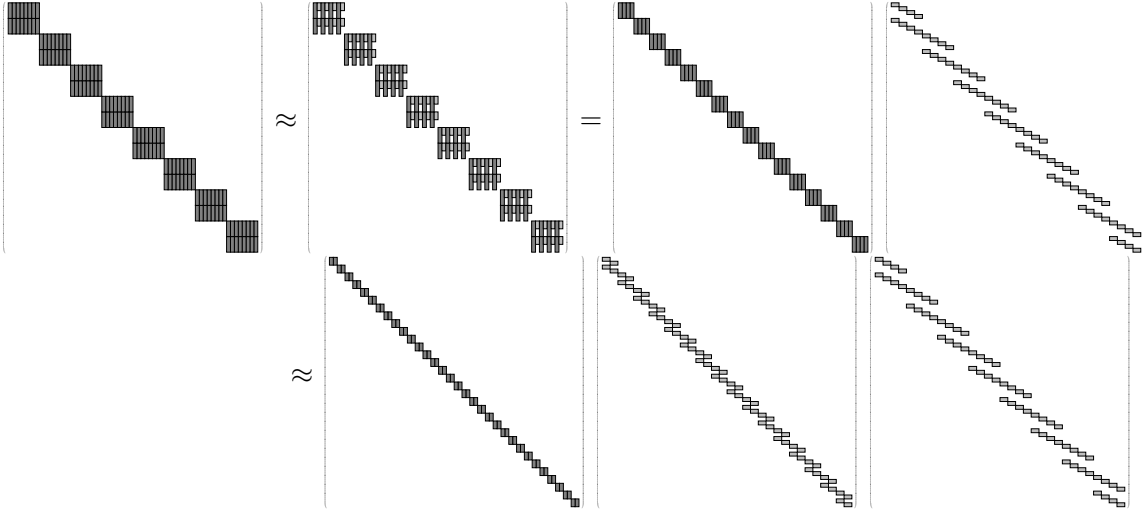
\begin{figure}[htp]
\begin{minipage}{\textwidth}
\centering
\resizebox{3.5cm}{!}{
\begin{tikzpicture}[baseline=-0.5ex]
      \tikzset{every left delimiter/.style={xshift=-1ex},every right delimiter/.style={xshift=1ex}}
      \matrix (mat) [matrix of math nodes, left delimiter=(, right delimiter=)] {
\draw;
\foreach \i in {0,0.5,...,3.5} {
    \foreach \j in {0,4,...,28} {
        \draw[fill=gray] (32-\i-\j,\j) rectangle (32-\i-\j-0.5,\j+2);
        \draw[fill=gray] (32-\i-\j,\j+2) rectangle (32-\i-\j-0.5,\j+4);
    }
}
\\
      };
\end{tikzpicture}
}
$\approx$
\resizebox{3.5cm}{!}{
\begin{tikzpicture}[baseline=-0.5ex]
      \tikzset{every left delimiter/.style={xshift=-1ex},every right delimiter/.style={xshift=1ex}}
      \matrix (mat) [matrix of math nodes, left delimiter=(, right delimiter=)] {
\draw;
\foreach \k in {0,4,...,28} {
\foreach \i in {0,...,3} {
    \foreach \j in {0,2} {
        \draw[fill=gray] (28-\k+\i,\j+\k) rectangle (28-\k+\i+0.5,\k+\j+2);
        \draw[fill=lightgray] (28-\k+\i+0.5,\j+\k+1) rectangle (28-\k+\i+1,\j+\k+2);
    }
}
}
\\
      };
\end{tikzpicture}
}
=
\resizebox{3.5cm}{!}{
\begin{tikzpicture}[baseline=-0.5ex]
      \tikzset{every left delimiter/.style={xshift=-1ex},every right delimiter/.style={xshift=1ex}}
      \matrix (mat) [matrix of math nodes, left delimiter=(, right delimiter=)] {
\draw;
\foreach \k in {0,2,...,30}{
\foreach \i in {0,0.5,...,1.5} {
    \draw[fill=gray] (32-\i-\k,\k) rectangle (32-\i-\k-0.5,\k+2);
}
}\\
      };
\end{tikzpicture}
}
\resizebox{3.5cm}{!}{
\begin{tikzpicture}[baseline=-0.5ex]
      \tikzset{every left delimiter/.style={xshift=-1ex},every right delimiter/.style={xshift=1ex}}
      \matrix (mat) [matrix of math nodes, left delimiter=(, right delimiter=)] {
\draw;
\foreach \k in {0,4,...,28}{
\foreach \i in {0,0.5,...,1.5} {
    \draw[fill=lightgray] (32-2*\i-\k,\k+\i) rectangle (32-2*\i-\k-1,\k+\i+0.5);
    \draw[fill=lightgray] (32-2*\i-\k,\k+\i+2) rectangle (32-2*\i-\k-1,\k+\i+2.5);
}
}\\
      };
\end{tikzpicture}
}\\
\resizebox{3.5cm}{!}{
}
$\approx$
\resizebox{3.5cm}{!}{
\begin{tikzpicture}[baseline=-0.5ex]
      \tikzset{every left delimiter/.style={xshift=-1ex},every right delimiter/.style={xshift=1ex}}
      \matrix (mat) [matrix of math nodes, left delimiter=(, right delimiter=)] {
\draw;
\foreach \k in {0,1,...,31}{
\foreach \i in {0,0.5} {
    \draw[fill=gray] (32-\i-\k,\k) rectangle (32-\i-\k-0.5,\k+1);
}
}\\
      };
\end{tikzpicture}
}
\resizebox{3.5cm}{!}{
\begin{tikzpicture}[baseline=-0.5ex]
      \tikzset{every left delimiter/.style={xshift=-1ex},every right delimiter/.style={xshift=1ex}}
      \matrix (mat) [matrix of math nodes, left delimiter=(, right delimiter=)] {
\draw;
\foreach \k in {0,2,...,30}{
\foreach \i in {0,0.5} {
    \draw[fill=lightgray] (32-2*\i-\k,\k+\i) rectangle (32-2*\i-\k-1,\k+\i+0.5);
    \draw[fill=lightgray] (32-2*\i-\k,\k+\i+1) rectangle (32-2*\i-\k-1,\k+\i+1.5);
}
}\\
      };
\end{tikzpicture}
}
\resizebox{3.5cm}{!}{
\begin{tikzpicture}[baseline=-0.5ex]
      \tikzset{every left delimiter/.style={xshift=-1ex},every right delimiter/.style={xshift=1ex}}
      \matrix (mat) [matrix of math nodes, left delimiter=(, right delimiter=)] {
\draw;
\foreach \k in {0,4,...,28}{
\foreach \i in {0,0.5,...,1.5} {
    \draw[fill=lightgray] (32-2*\i-\k,\k+\i) rectangle (32-2*\i-\k-1,\k+\i+0.5);
    \draw[fill=lightgray] (32-2*\i-\k,\k+\i+2) rectangle (32-2*\i-\k-1,\k+\i+2.5);
}
}\\
      };
\end{tikzpicture}
}\\
\end{minipage}
\caption{The recursive factorization of $U^3$
    in Figure~\ref{fig:compression-1D}.
    Gray factors are matrices inheriting the complementary low-rank
    property.
    Top: left matrix: $U^3$ with each diagonal block
    partitioned into smaller blocks according to
    Equation \eqref{eqn:half} as indicated by black rectangles;
    middle-left matrix: low-rank approximations of submatrices
    in $U^3$ given by Equation \eqref{eqn:RSh};
    middle right matrix: $U^4$;
    right matrix: $G^3$.
    Bottom: $U^4$ in the first row is further factorized into
    $U^4\approx U^5G^4$, giving $U^3 \approx U^5 G^4 G^3$.}
\label{fig:compression-1D2}
\end{figure}

Since there are $\O(1)$ nonzero entries in each $G^{h}_{i,j}$
and there are $\O(N)$ such submatrices,
there are only $\O(N)$ nonzero entries in $G^{h}$.
See Figure~\ref{fig:compression-1D2} top
for an example of the factorization $U^h\approx U^{h+1}G^h$
for the left factor $U^h$ with $L=6$,
$h=3$ and $r=1$ in Figure~\ref{fig:compression-1D}.

Similarly, for any $\ell$ between $h$ and $L-1$, we can factorize
$U^\ell\approx U^{\ell+1}G^{\ell}$, because the columns in $\left(
U^{\ell,t}_{i,2j} U^{\ell,t}_{i,2j+1} \right)$ and $\left(
U^{\ell,b}_{i,2j} U^{\ell,b}_{i,2j+1} \right)$ are in the column space
of the numerically low-rank matrices $K^{\ell+1}_{2i,j}$ and
$K^{\ell+1}_{2i+1,j}$, respectively. Computing the rank-$r$
approximations via the standard truncated {SVD} and transforming them
into the form of \eqref{eq:lowrankSVD3} give
\begin{equation}
\begin{pmatrix}
U^{\ell,t}_{i,2j} &  U^{\ell,t}_{i,2j+1}
\end{pmatrix}
\approx U^{\ell+1}_{2i,j} G^{\ell}_{2i,j}
\text{   and   }
\begin{pmatrix}
U^{\ell,b}_{i,2j} & U^{\ell,b}_{i,2j+1}
\end{pmatrix}
\approx U^{\ell+1}_{2i+1,j} G^{\ell}_{2i+1,j}
\label{eq:RSh2}
\end{equation}
for $i=0,1,\dots,2^\ell-1$ and $j=0,1,\dots,2^{L-\ell-1}-1$.  After assembling
these factorizations together, we obtain
\begin{equation*}
U^\ell\approx U^{\ell+1}G^{\ell}=
\begin{pmatrix}
U^{\ell+1}_{0} & & &\\
 & U^{\ell+1}_{1} & &\\
 & & \ddots &\\
 & & & U^{\ell+1}_{2^{\ell+1}-1}
\end{pmatrix}
\begin{pmatrix}
G^{\ell}_{0} & & &\\
 & G^{\ell}_{1} & &\\
 & & \ddots &\\
 & & & G^{\ell}_{2^\ell-1}
\end{pmatrix},
\end{equation*}
where
\begin{equation*}
U_{i}^{\ell+1}=
\begin{pmatrix}
U^{\ell+1}_{i,0} & U^{\ell+1}_{i,1} & \cdots & U^{\ell+1}_{i,2^{L-\ell-1}-1}
\end{pmatrix}
\end{equation*}
for $i=0,1,\dots,2^{\ell+1}-1$, and
\begin{equation*}
G^{\ell}_{i}=
\begin{pmatrix}
G^{\ell}_{2i,0} & & &\\
 & G^{\ell}_{2i,1} & &\\
 & & \ddots &\\
 & & & G^{\ell}_{2i,2^{L-\ell-1}-1}\\
\midrule
G^{\ell}_{2i+1,0} & & &\\
 & G^{\ell}_{2i+1,1} & &\\
 & & \ddots &\\
 & & & G^{\ell}_{2i+1,2^{L-\ell-1}-1}
\end{pmatrix}
\end{equation*}
for $i=0,1,\dots,2^\ell-1$.

After $L-h$ steps of recursive factorizations
\begin{equation*}
U^\ell \approx U^{\ell+1}G^{\ell}
\end{equation*}
for $\ell=h,h+1,\dots,L-1$,
we obtain the recursive factorization of $U^h$ as
\begin{equation}
\label{eqn:RFUh}
U^h\approx U^{L}G^{L-1}\cdots G^{h}.
\end{equation}
See Figure~\ref{fig:compression-1D2} bottom for an example
of a recursive factorization
for the left factor $U^h$ with $L=6$, $h=3$ and $r=1$
in Figure~\ref{fig:compression-1D}.

Similar to the analysis of $G^{h}$, it is also easy to check that
there are only $\O(N)$ nonzero entries in each $G^{\ell}$
in \eqref{eqn:RFUh}.  Since there are $\O(N)$ diagonal blocks in
$U^L$ and each block contains $\O(1)$ entries, there is
$\O(N)$ nonzero entries in $U^L$.

\subsubsection{Recursive factorization of $V^h$}
\label{sec:RCV}

\begin{figure}[htp]
\begin{minipage}{\textwidth}
\centering
\resizebox{3.5cm}{!}{
\begin{tikzpicture}[baseline=-0.5ex]
      \tikzset{every left delimiter/.style={xshift=-1ex},every right delimiter/.style={xshift=1ex}}
      \matrix (mat) [matrix of math nodes, left delimiter=(, right delimiter=)] {
      \draw;
\foreach \i in {0,0.5,...,3.5} {
    \foreach \j in {0,4,...,28} {
        \draw[fill=gray] (\j,32-\i-\j) rectangle (\j+2,32-\i-\j-0.5);
        \draw[fill=gray] (\j+2,32-\i-\j) rectangle (\j+4,32-\i-\j-0.5);
    }
}
\\
      };
\end{tikzpicture}
}
=
\resizebox{3.5cm}{!}{
\begin{tikzpicture}[baseline=-0.5ex]
      \tikzset{every left delimiter/.style={xshift=-1ex},every right delimiter/.style={xshift=1ex}}
      \matrix (mat) [matrix of math nodes, left delimiter=(, right delimiter=)] {
\draw;
\foreach \k in {0,4,...,28}{
\foreach \i in {0,0.5,...,1.5} {
    \draw[fill=lightgray] (\k+\i,32-2*\i-\k) rectangle (\k+\i+0.5,32-2*\i-\k-1);
    \draw[fill=lightgray] (\k+\i+2,32-2*\i-\k) rectangle (\k+\i+2.5,32-2*\i-\k-1);
}
}\\
      };
\end{tikzpicture}
}
\resizebox{3.5cm}{!}{
\begin{tikzpicture}[baseline=-0.5ex]
      \tikzset{every left delimiter/.style={xshift=-1ex},every right delimiter/.style={xshift=1ex}}
      \matrix (mat) [matrix of math nodes, left delimiter=(, right delimiter=)] {
\draw;
\foreach \k in {0,2,...,30}{
\foreach \i in {0,0.5} {
    \draw[fill=lightgray] (\k+\i,32-2*\i-\k) rectangle (\k+\i+0.5,32-2*\i-\k-1);
    \draw[fill=lightgray] (\k+\i+1,32-2*\i-\k) rectangle (\k+\i+1.5,32-2*\i-\k-1);
}
}\\
      };
\end{tikzpicture}
}
\resizebox{3.5cm}{!}{
\begin{tikzpicture}[baseline=-0.5ex]
      \tikzset{every left delimiter/.style={xshift=-1ex},every right delimiter/.style={xshift=1ex}}
      \matrix (mat) [matrix of math nodes, left delimiter=(, right delimiter=)] {
\draw;
\foreach \k in {0,1,...,31}{
\foreach \i in {0,0.5} {
    \draw[fill=gray] (\k,32-\i-\k) rectangle (\k+1,32-\i-\k-0.5);
}
}\\
      };
\end{tikzpicture}
}\\
\end{minipage}
\caption{The recursive factorization $(V^3)^* \approx (H^3)^*(H^4)^*(V^5)^*$ of $(V^3)^*$
    in Figure~\ref{fig:compression-1D}.}
\label{fig:compression-1D3}
\end{figure}
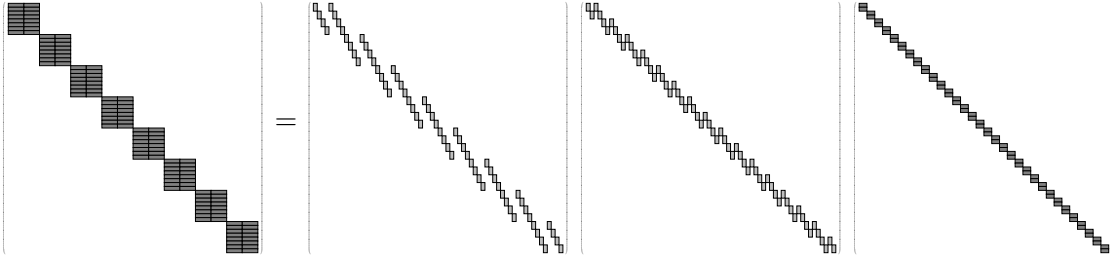

The recursive factorization of $V^h$ is similar to the one of $U^h$.
In each step of the factorization
\[
(V^\ell)^*\approx (H^{\ell})^* (V^{\ell+1})^*,
\]
we take advantage of the low-rank property of the row space of
$K^{L-\ell-1}_{i,2j}$ and $K^{L-\ell-1}_{i,2j+1}$ to obtain rank-$r$
approximations. Applying the exact same procedure of Section
\ref{sec:RCU} now to $V^\ell$ leads to the recursive factorization
$V^h \approx V^L H^{L-1} \cdots H^{h}$, or equivalently
\begin{equation}
\label{eqn:RFVh}
(V^h)^* \approx (H^{h})^* \cdots (H^{L-1})^* (V^L)^*,
\end{equation}
with all factors containing only $\O(N)$ nonzero entries.  See
Figure~\ref{fig:compression-1D3} for an example of a recursive
factorization $(V^h)^* \approx (H^{h})^* \cdots (H^{L-2})^*
(V^{L-1})^*$ for the left factor $V^h$ with $L=6$, $h=3$ and $r=1$ in
Figure~\ref{fig:compression-1D}.

Given the recursive factorization of $U^h$ and $(V^h)^*$ in
\eqref{eqn:RFUh} and \eqref{eqn:RFVh}, we reach the butterfly
factorization
\begin{equation}
K \approx U^{L}G^{L-1}\cdots G^{h} M^h (H^{h})^* \cdots (H^{L-1})^* (V^L)^*,
\end{equation}
where all factors are sparse matrices with $\O(N)$ nonzero entries.
For a given input vector $g\in\bbC^N$, the $\O(N^2)$ matrix-vector
multiplication $u = Kg$ can be approximated by a sequence of $\O(\log
N)$ sparse matrix-vector multiplications given by the butterfly
factorization.

\subsection{Complexity analysis}
\label{sec:Complexity}
The complexity analysis of the construction of a butterfly
factorization naturally consists of two parts: the middle level
factorization and the recursive factorization.

The complexity of the middle level factorization depends on which one
of the cases is under consideration.
\begin{itemize}
\item For the first case, the approximate SVDs are determined by the
  application of $K$ and $K^*$ to Gaussian random matrices in
  $\bbC^{N\times N^{1/2}(r+p)}$ and the rank-$r$ approximations of
  $K^h_{ij}$ for each $(i,j)$ pair.  Assume that each matrix-vector
  multiplication by $K$ or $K^*$ via the given black-box routines
  requires $\O(C_K(N))$ operations (which is at least $\O(N)$).  Then
  the dominant cost is due to applying $K$ and $K^*$ $\O(N^{1/2})$
  times, which yields an overall computational complexity of
  $\O(C_K(N)N^{1/2})$.
\item In the second case, the approximate SVDs are computed via random
  sampling for each $K^h_{ij}$ of the $\O(N)$ pairs $(i,j)$. The
  complexity of performing randomized sampling for each such block is
  $\O(N^{1/2})$.  Hence, the overall computational complexity is
  $\O(N^{3/2})$.
\end{itemize}

In the recursive factorization, $U^{\ell}$ at level $\ell$ consists of
$\O(2^{\ell})$ diagonal blocks of size $\O(N/2^{\ell})\times \O(
N/2^{\ell})$.  In each diagonal block, there are $\O(N/2^{\ell})$
factorizations in \eqref{eq:RSh2}.  Since the operation complexity of
performing one factorization in \eqref{eq:RSh2} is $\O (N /
2^{\ell})$, it takes $\O(N^{2} / 2^{\ell})$ operations to factorize
$U^{\ell}$.  Summing up the operations at all levels gives the total
complexity for recursively factorizing $U^h$:
\begin{equation}\label{eq:complexityU}
\sum^{L-1}_{\ell=h}\O(N^{2}/2^{\ell})
= \O(N^{3/2}).
\end{equation}
Similarly, the operation complexity for recursively compressing $V^h$
is also $\O(N^{3/2})$.

The memory peak of the butterfly factorization occurs in the middle
level factorization since we have to store the initial factorization
in \eqref{eq:KeqUMV}.  There are $\O(N^{3/2})$ nonzero entries in
$U^h$ and $V^h$, and $\O(N)$ in $M^h$.  Hence, the total memory
complexity is $\O(N^{3/2})$. The total operation complexity for
constructing the butterfly factorization is summarized in Table
\ref{tab:complexity}.

\begin{table}[htb]
\centering
\begin{tabular}{llcc}
\toprule
 & & \parbox{2cm}{\centering Randomized\\{SVD}}
 & \parbox{2cm}{\centering Randomized\\sampling} \\
\toprule
\multirow{5}{3cm}{Factorization Complexity} &
\parbox{2.5cm}{Middle level\\factorization}
& $\O(C_K(N)N^{1/2})$ & $\O(N^{3/2})$ \\
\cmidrule(r){2-4}
& \parbox{2.5cm}{Recursive\\factorization} & \multicolumn{2}{c}{$\O(N^{3/2})$}\\
\cmidrule(r){2-4}
& Total & $\O(C_K(N)N^{1/2})$ & $\O(N^{3/2})$\\
\toprule
\parbox{3cm}{Memory\\Complexity} & & $\O(N^{3/2})$ & $\O(N\log N)$ \\
\toprule
\parbox{3cm}{Application\\Complexity} & & \multicolumn{2}{c}{$\O(N\log N)$}\\
\bottomrule
\end{tabular}
\caption{Computational complexity and memory complexity of the
  butterfly factorization. $C_K(N)$ is the operation complexity of one
  application of $K$ or $K^*$. In most of the cases encountered,
  $C_K(N)=O(N\log N)$.}
\label{tab:complexity}
\end{table}

It is worth pointing out that the memory complexity can be reduced to
$\O(N\log N)$, when we apply the randomized sampling method to
construct each block in the initial factorization in \eqref{eq:KeqUMV}
separately.  Instead of factorizing $U^h$ and $V^h$ at the end of the
middle level factorization, we can factorize the left and right
factors $U^h_i$ and $V^h_i$ in \eqref{eq:KeqUMV} on the fly to avoid
storing all factors in \eqref{eq:KeqUMV}.  For a fixed $i$, we
generate $U^h_{i}$ from $K^h_{ij}$ for all $j$, and recursively
factorize $U^h_{i}$.  The memory cost is $\O(N)$ for storing $U^h_{i}$
and $\O(N^{1/2}\log N)$ for storing the sparse matrices after its
recursive factorization.  Repeating this process for $i = 1, \dots,
N^{1/2}$ gives the complete factorization of $U^h$.  The factorization
of $V^h$ is conducted similarly.  The total memory complexity is
$\O(N\log N)$.

The operation and memory complexity for the application of
the butterfly factorization are governed
by the number of nonzero entries
in the factorization: $\O(N\log N)$.

\section{Numerical results}
\label{sec:results}
This section presents three numerical examples to demonstrate the
effectiveness of the algorithms proposed above.  The first example is
an FIO in \cite{FIO09} and the second example is a special function
transform in \cite{Butterfly2}. Both examples provide an explicit
kernel function that becomes a one-dimensional complementary low-rank
matrix after discretization.  This allows us to apply the butterfly
factorization construction algorithm with random sampling. The
computational complexity and the memory cost are $\O(N^{3/2})$ and
$\O(N\log N)$ in this case.

The third example is a composition of two FIOs for which an explicit
kernel function of their composition is not available.  Since we can
apply either the butterfly algorithm in \cite{FIO09} or the butterfly
factorization to evaluate these FIOs one by one, a fast algorithm for
computing the composition is available. We apply the butterfly
factorization construction algorithm with random matrix-vector
multiplication to this example which requires $\O(N^{3/2}\log N)$
operations and $\O(N^{3/2})$ memory cost.

Our implementation is in MATLAB.  The numerical results were obtained
on a server computer with a 2.0 {GHz CPU}. The additive oversampling
parameter is $p=5$ and the multiplicative oversampling parameter is
$q=3$.

Let $\{u^d(x),x\in X\}$ and $\{u^a(x),x\in X\}$ denote the results
given by the direct matrix-vector multiplication and the butterfly
factorization.  The accuracy of applying the butterfly factorization
algorithm is estimated by the following relative error
\begin{equation}
\epsilon^a = \sqrt{\cfrac{\sum_{x\in S}|u^a(x)-u^d(x)|^2}
{\sum_{x\in S}|u^d(x)|^2}},
\end{equation}
where $S$ is a point set of size $256$ randomly sampled from $X$.

\paragraph{Example 1.}
Our first example is to evaluate a one-dimensional FIO of the
following form:
\begin{equation}\label{eq:example1}
u(x) = \int_{\mathbb{R}}e^{2\pi \imath \Phi(x,\xi)}\widehat{f}(\xi)d\xi,
\end{equation}
where $\widehat{f}$ is the {Fourier} transform of $f$,
and $\Phi(x,\xi)$ is a phase function given by
\begin{equation}
\Phi(x,\xi) = x\cdot \xi + c(x)|\xi|,~~~c(x) = (2+\sin(2\pi x))/8.
\label{eqn:1D-FIO-kernal}
\end{equation}
The discretization of \eqref{eq:example1} is
\begin{equation}
\label{eqn:1D-FIO}
u(x_i) = \sum_{\xi_j}e^{2\pi \imath \Phi(x_i,\xi_j)}\widehat{f}(\xi_j),
\quad i,j=1,2,\dots,N,
\end{equation}
where $\{x_i\}$ and $\{\xi_j\}$ are uniformly distributed points
in $[0,1)$ and $[-N/2,N/2)$ following
\begin{equation}\label{eqn:1D-xandxi}
x_i = (i-1)/N \text{ and } \xi_j = j-1-N/2.
\end{equation}
\eqref{eqn:1D-FIO} can be represented in a matrix form as $u=Kg$,
where $u_i=u(x_i)$, $K_{ij} = e^{2\pi\imath \Phi(x_i,\xi_j)}$ and $g_j
= \widehat{f}(\xi_j)$.  The matrix $K$ satisfies the complementary
low-rank property as proved in \cite{FIO09,FIO14}.  The explicit
kernel function of $K$ allows us to use the construction algorithm
with random sampling. Table~\ref{tab:1D-fio} summarizes the results of
this example for different grid sizes $N$ and truncation ranks $r$.
\begin{table}[htp]
\centering
\begin{tabular}{rccccc}
\toprule
   $N,r$ & $\epsilon^a$&$T_{Factor}(min)$&$T_d(sec)$&$T_a(sec)$
                                       &$T_d/T_a$\\
\toprule
  1024,4 & 2.49e-05 & 2.92e-01 & 2.30e-01 & 3.01e-02 & 7.65e+00 \\
  4096,4 & 4.69e-05 & 1.62e+00 & 2.64e+00 & 4.16e-02 & 6.35e+01 \\
 16384,4 & 5.77e-05 & 1.22e+01 & 2.28e+01 & 1.84e-01 & 1.24e+02 \\
 65536,4 & 6.46e-05 & 8.10e+01 & 2.16e+02 & 1.02e+00 & 2.12e+02 \\
262144,4 & 7.13e-05 & 4.24e+02 & 3.34e+03 & 4.75e+00 & 7.04e+02 \\
\toprule
  1024,6 & 1.57e-08 & 1.81e-01 & 1.84e-01 & 1.20e-02 & 1.54e+01 \\
  4096,6 & 3.64e-08 & 1.55e+00 & 2.56e+00 & 6.42e-02 & 3.98e+01 \\
 16384,6 & 6.40e-08 & 1.25e+01 & 2.43e+01 & 3.01e-01 & 8.08e+01 \\
 65536,6 & 6.53e-08 & 9.04e+01 & 2.04e+02 & 1.77e+00 & 1.15e+02 \\
262144,6 & 6.85e-08 & 5.45e+02 & 3.68e+03 & 8.62e+00 & 4.27e+02 \\
\toprule
  1024,8 & 5.48e-12 & 1.83e-01 & 1.78e-01 & 1.63e-02 & 1.09e+01 \\
  4096,8 & 1.05e-11 & 1.98e+00 & 2.71e+00 & 8.72e-02 & 3.11e+01 \\
 16384,8 & 2.09e-11 & 1.41e+01 & 3.34e+01 & 5.28e-01 & 6.33e+01 \\
 65536,8 & 2.62e-11 & 1.17e+02 & 2.10e+02 & 2.71e+00 & 7.75e+01 \\
262144,8 & 4.13e-11 & 6.50e+02 & 3.67e+03 & 1.52e+01 & 2.42e+02 \\
\bottomrule
\end{tabular}
\caption{Numerical results for the FIO given in \eqref{eqn:1D-FIO}.
  $N$ is the size of the matrix; $r$ is the fixed rank in the low-rank
  approximations; $T_{Factor}$ is the factorization time of the
  butterfly factorization; $T_{d}$ is the running time of the direct
  evaluation; $T_{a}$ is the application time of the butterfly
  factorization; $T_d/T_a$ is the speedup factor.  }
\label{tab:1D-fio}
\end{table}

\paragraph{Example 2.}
Next, we provide an example of a special function transform.
This example can be further applied to accelerate the Fourier-Bessel transform
that is important in many real applications.
Following the standard notation, we denote the {Hankel} function of
the first kind of order $m$ by $H^{(1)}_m$.  When $m$ is an integer,
$H^{(1)}_m$ has a singularity at the origin and a branch cut along the negative
real axis.  We are interested in evaluating the sum of {Hankel}
functions over different orders,
\begin{equation}
  \label{eqn:1D-special}
  u(x_i) = \sum_{j=1}^{N}H^{(1)}_{j-1}(x_i)g_j,
  \quad i=1,2,\dots,N,
\end{equation}
which is analogous to expansion in orthogonal polynomials.  The points
$x_i$ are defined via the formula,
\begin{equation}
  x_i = N+\cfrac{2\pi}{3}(i-1)
\end{equation}
which are bounded away from zero. It is demonstrated in
\cite{Butterfly2} that \eqref{eqn:1D-special} can be represented via
$u=Kg$ where $K$ satisfies the complementary low-rank property,
$u_i=u(x_i)$ and $K_{ij} = H^{(1)}_{j-1}(x_i)$.  The entries of matrix
$K$ can be calculated efficiently and the construction
algorithm with random sampling is applied to accelerate the evaluation
of the sum \eqref{eqn:1D-special}.  Table~\ref{tab:1D-special}
summarizes the results of this example for different grid sizes $N$
and truncation ranks $r$.

\begin{table}[htp]
  \centering
  \begin{tabular}{rccccc}
    \toprule
    $N,r$ & $\epsilon^a$&$T_{Factor}(min)$&$T_d(sec)$&$T_a(sec)$
    &$T_d/T_a$ \\
\toprule
  1024,4 & 2.35e-06 & 8.78e-01 & 8.30e-01 & 1.06e-02 & 7.86e+01 \\
  4096,4 & 5.66e-06 & 5.02e+00 & 5.30e+00 & 2.83e-02 & 1.87e+02 \\
 16384,4 & 6.86e-06 & 3.04e+01 & 5.51e+01 & 1.16e-01 & 4.76e+02 \\
 65536,4 & 7.04e-06 & 2.01e+02 & 7.59e+02 & 6.38e-01 & 1.19e+03 \\
\toprule
  1024,6 & 2.02e-08 & 4.31e-01 & 7.99e-01 & 9.69e-03 & 8.25e+01 \\
  4096,6 & 4.47e-08 & 6.61e+00 & 5.41e+00 & 4.52e-02 & 1.20e+02 \\
 16384,6 & 5.95e-08 & 4.19e+01 & 5.62e+01 & 1.61e-01 & 3.48e+02 \\
 65536,6 & 7.86e-08 & 2.76e+02 & 7.60e+02 & 1.01e+00 & 7.49e+02 \\
    \bottomrule
  \end{tabular}
  \caption{Numerical results with the matrix given by
    \eqref{eqn:1D-special}.}
  \label{tab:1D-special}
\end{table}

From Table \ref{tab:1D-fio} and \ref{tab:1D-special}, we note that the
accuracy of the butterfly factorization is well controlled by the max
rank $r$.  For a fixed rank $r$, the accuracy is almost independent of
$N$.  In practical applications, one can set the desired $\epsilon$
ahead and increase the truncation rank $r$ until the relative error
reaches $\epsilon$.

The tables for Example 1 and Example 2 also provide
numerical evidence for the asymptotic
complexity of the proposed algorithms. The construction algorithm
based on random sampling is of computational complexity $\O(N^{3/2})$.
When we quadruple the problem size, the running time of the
construction sextuples and is better than we expect.  The reason is
that in the random sampling method, the computation of a middle matrix
requires pseudo-inverses of $r \times r$ matrices whose complexity is
$\O(r^3)$ with a large prefactor.  Hence, when $N$ is not large, the
running time will be dominated by the $\O(r^3N)$ computation of middle
matrices.  The numbers also show that the application
complexity of the butterfly factorization is $\O(N\log N)$
with a prefactor much smaller than
the butterfly algorithm with {Chebyshev} interpolation \cite{FIO09}.
In example 1, when the relative error is $\epsilon\approx 10^{-5}$,
the butterfly factorization truncates the low-rank submatrices with rank 4
whereas the butterfly algorithm with {Chebyshev} interpolation
uses 9 {Chebyshev} grid points.
The speedup factors are 200 on average.

\paragraph{Example 3.}
In this example, we consider a composition of two {FIOs}, which is the
discretization of the following operator
\begin{equation}
  \label{eqn:composition}
  u(x) = \int_{\mathbb{R}}e^{2\pi \imath \Phi_2(x,\eta)}\int_{\mathbb{R}}
  e^{-2\pi \imath y\eta}\int_{\mathbb{R}}e^{2\pi \imath \Phi_1(y,\xi)}
  \widehat{f}(\xi)d\xi dy d\eta.
\end{equation}
For simplicity, we consider the same phase function
$\Phi_1=\Phi_2=\Phi$ as given by \eqref{eqn:1D-FIO-kernal}.  By the
discussion of Example $1$ for one FIO, we know the discrete analog of
the composition \eqref{eqn:composition} can be represented as
\begin{equation*}
  u=KFKF f =: KFK g,\quad \text{with}\; g=Ff,
\end{equation*}
where $F$ is the standard Fourier transform in matrix form, $K$ is the
same matrix as in Example $1$, $u_i=u(x_i)$, and $g_j =
\widehat{f}(\xi_j)$.  Under mild assumptions as discussed in
\cite{Theory}, the composition of two {FIOs} is an {FIO}.  Hence, the
new kernel matrix $\tilde{K} = K F K$ again satisfies the
complementary low-rank property, though typically with slightly
increased ranks.

Notice that it is not reasonable to compute the matrix $\tilde{K}$
directly.  However, we have the fast Fourier transform (FFT) to apply
$F$ and the butterfly factorization that we have built for $K$ in
Example $1$ to apply $K$. Therefore, the construction algorithm with
random matrix-vector multiplication is applied to factorize
$\tilde{K}$.

Since the direct evaluation of each $u_i$ takes $\O(N^2)$ operations,
the exact solution $\{u^d_i\}_{i\in S}$ for a selected set $S$ is
infeasible for large $N$.  We apply the butterfly factorization of $K$
and the FFT to evaluate $\{u_i\}_{i\in S}$ as an approximation to the
exact solution $\{u^d_i\}_{i\in S}$.  These approximations are
compared to the results $\{u^a_i\}_{i\in S}$ that are given by
applying the butterfly factorization of $\tilde{K}$.
Table~\ref{tab:1D-mfio}
summarizes the results of this example for different grid sizes $N$
and truncation ranks $r$.
\begin{table}[htp]
\centering
\begin{tabular}{rccccc}
  \toprule
  $N,r$ & $\epsilon^a$&$T_{Factor}(min)$&$T_d(sec)$&$T_a(sec)$
  &$T_d/T_a$ \\
\toprule
  1024,4 & 1.40e-02 & 3.26e-01 & 3.64e-01 & 4.74e-03 & 7.69e+01 \\
  4096,4 & 1.96e-02 & 4.20e+00 & 6.59e+00 & 2.52e-02 & 2.62e+02 \\
 16384,4 & 2.34e-02 & 4.65e+01 & 3.75e+01 & 1.15e-01 & 3.25e+02 \\
 65536,4 & 2.18e-02 & 4.33e+02 & 3.73e+02 & 6.79e-01 & 5.49e+02 \\
\toprule
  1024,8 & 6.62e-05 & 3.65e-01 & 3.64e-01 & 8.25e-03 & 4.42e+01 \\
  4096,8 & 8.67e-05 & 4.94e+00 & 6.59e+00 & 5.99e-02 & 1.10e+02 \\
 16384,8 & 1.43e-04 & 6.23e+01 & 3.75e+01 & 3.47e-01 & 1.08e+02 \\
 65536,8 & 1.51e-04 & 6.91e+02 & 3.73e+02 & 1.76e+00 & 2.12e+02 \\
\toprule
  1024,12 & 1.64e-08 & 4.79e-01 & 3.64e-01 & 1.48e-02 & 2.46e+01 \\
  4096,12 & 1.05e-07 & 6.35e+00 & 6.59e+00 & 1.12e-01 & 5.88e+01 \\
 16384,12 & 2.55e-07 & 7.58e+01 & 3.75e+01 & 7.64e-01 & 4.91e+01 \\
 65536,12 & 2.69e-07 & 7.63e+02 & 3.73e+02 & 4.39e+00 & 8.49e+01 \\
\bottomrule
\end{tabular}
\caption{Numerical results for the composition of two FIOs.}
\label{tab:1D-mfio}
\end{table}

Table \ref{tab:1D-mfio} shows the numerical results of the butterfly
factorization of $\tilde{K}$.  The accuracy improves as we increase the
truncation rank $r$.  Comparing Table \ref{tab:1D-mfio} with Table
\ref{tab:1D-fio}, we notice that, for a fixed accuracy, the rank used
in the butterfly factorization of the composition of {FIOs} should be
larger than the rank used in a single {FIO} butterfly factorization.
This is expected since the composition is in general more complicated
than the individual FIOs. $T_{Factor}$ grows on average by a factor of
ten when we quadruple the problem size. This agrees with the estimated
$\O(N^{3/2}\log N)$ computational complexity for constructing the
butterfly factorization.  The column $T_a$ shows that the empirical
application time of our factorization is close to the estimated
complexity $O(N\log N)$.

\section{Conclusion and discussion}
\label{sec:conclusion}
This paper introduces a butterfly factorization as a data-sparse
approximation of complementary low-rank matrices.  More precisely, it
represents such an $N\times N$ dense matrix as a product of $\O(\log
N)$ sparse matrices. The factorization can be built efficiently if
either a fast algorithm for applying the matrix and its adjoint is
available or an explicit expression for the entries of the matrix is
given. The butterfly factorization gives rise to highly efficient
matrix-vector multiplications with $\O(N\log N)$ operation and memory
complexity. The butterfly factorization is also useful when an
existing butterfly algorithm is repeatedly applied, because the
application of the butterfly factorization is significantly faster
than pre-existing butterfly algorithms.

The method proposed here is a first step in computing data-sparse
approximations for butterfly algorithms. As we mentioned earlier, the
proposed butterfly factorization can be easily generalized to higher
dimensions, which is especially relevant in imaging science. Another
interesting direction is to invert a matrix via the butterfly
factorization. While the numerical results of this paper include a
couple of examples, it is natural to consider other important class of
transforms, such as the non-uniform Fourier transform and the Legendre
functions associated with the spherical harmonic transform.

{\bf Acknowledgments.}  Y. Li, H. Yang and L. Ying were partially
supported by the National Science Foundation under award DMS-1328230
and the U.S. Department of Energy's Advanced Scientific Computing
Research program under award DE-FC02-13ER26134/DE-SC0009409.
E. Martin was supported in part by DOE grant number DE-FG02-97ER25308.
K. Ho was supported by the National Science Foundation under award
DMS-1203554.

\bibliographystyle{abbrv} \bibliography{ref}

\begin{thebibliography}{10}

\bibitem{FIO09}
E.~J. Cand{\`e}s, L.~Demanet, and L.~Ying.
\newblock A fast butterfly algorithm for the computation of {Fourier} integral
  operators.
\newblock {\em Multiscale Modeling and Simulation}, 7(4):1727--1750, 2009.

\bibitem{FIO12}
L.~Demanet and L.~Ying.
\newblock Fast wave computation via {Fourier} integral operators.
\newblock {\em Mathematics of Computation}, 81:1455--1486, 2012.

\bibitem{Butterfly4}
B.~Engquist and L.~Ying.
\newblock A fast directional algorithm for high frequency acoustic scattering
  in two dimensions.
\newblock {\em Communications in Mathematical Sciences}, 7(2):327--345, 06
  2009.

\bibitem{HMat}
W.~Hackbusch.
\newblock {A sparse matrix arithmetic based on $\mathcal{H}$-matrices. I.
  Introduction to $\mathcal{H}$-matrices}.
\newblock {\em Computing}, 62(2):89--108, 1999.

\bibitem{HSS2}
W.~Hackbusch and S.~B\"{o}rm.
\newblock Data-sparse approximation by adaptive {$\mathcal{H}^2$}-matrices.
\newblock {\em Computing}, 69(1):1--35, 2002.

\bibitem{Rec1}
N.~Halko, P.~Martinsson, and J.~Tropp.
\newblock Finding structure with randomness: Probabilistic algorithms for
  constructing approximate matrix decompositions.
\newblock {\em SIAM Review}, 53(2):217--288, 2011.

\bibitem{Theory}
L.~H\"{o}rmander.
\newblock {Fourier} integral operators. {I}.
\newblock {\em Acta Mathematica}, 127(1):79--183, 1971.

\bibitem{FIO14}
Y.~Li, H.~Yang, and L.~Ying.
\newblock A multiscale butterfly aglorithm for {Fourier} integral operators.
\newblock {\em Multiscale Modeling and Simulation}, to appear.

\bibitem{Rec2}
E.~Liberty, F.~Woolfe, P.-G. Martinsson, V.~Rokhlin, and M.~Tygert.
\newblock Randomized algorithms for the low-rank approximation of matrices.
\newblock {\em Proc. Natl. Acad. Sci. USA}, 104(51):20167--20172, 2007.

\bibitem{HMatrix}
L.~Lin, J.~Lu, and L.~Ying.
\newblock Fast construction of hierarchical matrix representation from
  matrix-vector multiplication.
\newblock {\em J. Comput. Phys.}, 230(10):4071--4087, 2011.

\bibitem{mahoney2009cur}
M.~W. Mahoney and P.~Drineas.
\newblock {CUR} matrix decompositions for improved data analysis.
\newblock {\em Proceedings of the National Academy of Sciences},
  106(3):697--702, 2009.

\bibitem{HSSMatrix}
P.~G. Martinsson.
\newblock A fast randomized algorithm for computing a hierarchically
  semiseparable representation of a matrix.
\newblock {\em SIAM J. Matrix Anal. Appl.}, 32(4):1251--1274, 2011.

\bibitem{Butterfly1}
E.~Michielssen and A.~Boag.
\newblock A multilevel matrix decomposition algorithm for analyzing scattering
  from large structures.
\newblock {\em Antennas and Propagation, IEEE Transactions on},
  44(8):1086--1093, Aug 1996.

\bibitem{Butterfly2}
M.~O'Neil, F.~Woolfe, and V.~Rokhlin.
\newblock An algorithm for the rapid evaluation of special function transforms.
\newblock {\em Appl. Comput. Harmon. Anal.}, 28(2):203--226, 2010.

\bibitem{FIO13}
J.~Poulson, L.~Demanet, N.~Maxwell, and L.~Ying.
\newblock A parallel butterfly algorithm.
\newblock {\em SIAM J. Sci. Comput.}, 36(1):C49--C65, 2014.

\bibitem{wavemoth}
D.~S. Seljebotn.
\newblock Wavemoth-fast spherical harmonic transforms by butterfly matrix
  compression.
\newblock {\em The Astrophysical Journal Supplement Series}, 199(1):5, 2012.

\bibitem{InvRadon}
D.~O. Trad, T.~J. Ulrych, and M.~D. Sacchi.
\newblock Accurate interpolation with high-resolution time-variant {Radon}
  transforms.
\newblock {\em Geophysics}, 67(2):644--656, 2002.

\bibitem{Special2}
M.~Tygert.
\newblock Fast algorithms for spherical harmonic expansions, {III}.
\newblock {\em Journal of Computational Physics}, 229(18):6181 -- 6192, 2010.

\bibitem{Rec3}
F.~Woolfe, E.~Liberty, V.~Rokhlin, and M.~Tygert.
\newblock A fast randomized algorithm for the approximation of matrices.
\newblock {\em Applied and Computational Harmonic Analysis}, 25(3):335 -- 366,
  2008.

\bibitem{HSS1}
J.~Xia, S.~Chandrasekaran, M.~Gu, and X.~S. Li.
\newblock Fast algorithms for hierarchically semiseparable matrices.
\newblock {\em Numerical Linear Algebra with Applications}, 17(6):953--976,
  2010.

\bibitem{symbol3}
H.~Yang and L.~Ying.
\newblock A fast algorithm for multilinear operators.
\newblock {\em Applied and Computational Harmonic Analysis}, 33(1):148 -- 158,
  2012.

\bibitem{InvFIO}
B.~Yazici, L.~Wang, and K.~Duman.
\newblock Synthetic aperture inversion with sparsity constraints.
\newblock In {\em Electromagnetics in Advanced Applications (ICEAA), 2011
  International Conference on}, pages 1404--1407, Sept 2011.

\bibitem{Butterfly5}
L.~Ying.
\newblock Sparse {Fourier} transform via butterfly algorithm.
\newblock {\em SIAM J. Sci. Comput.}, 31(3):1678--1694, Feb. 2009.

\end{thebibliography}

\end{document}